\documentclass[12pt]{article}

\usepackage{amsmath,amsfonts,amssymb}

\newenvironment{demo}{\noindent {\sl Proof}. \ }{\qed \bigskip}
\newenvironment{expl}{\noindent {\bf Example} \ }{}
\newtheorem{stheorem}{Theorem }[subsection]
\newtheorem{sdefin}[stheorem]{Definition }
\newtheorem{sprop}[stheorem]{Proposition }
\newtheorem{slemma}[stheorem]{Lemma }
\newtheorem{scorol}[stheorem]{Corollary }

\numberwithin{equation}{section}

\def\half{\frac{1}{2}}
\def\noi{\noindent}
\def\R{\mathbb{R}}
\def\N{\mathbb{N}}
\def\ep{\varepsilon}

\def\omb{\overline \omega}
\def\al{\alpha}
\def\la{\lambda}
\def\si{\sigma}
\def\cB{{\cal B}}
\def\cF{{\cal F}}
\def\tix{\tilde x}
\def\esp{{\rm I\mskip-4mu E}}
\def\Hess{{\rm Hess \ }}
\def\hY{\hat Y}
\def\hZ{\hat Z}
\def\dX{\delta X}
\def\dY{\delta Y}
\def\dZ{\delta Z}
\def\qed{\hbox{$\vcenter{\vbox{
   \hrule height 0.4pt\hbox{\vrule width 0.4pt height 6pt
    \kern5pt\vrule width 0.4pt}\hrule height 0.4pt}}$}}

\begin{document}

\title{
\bigskip \bigskip
        Backward Stochastic Differential Equations on Manifolds}
\author{
\begin{tabular}{c}
Fabrice Blache \\
{\it{Laboratoire de Math\'ematiques}}                 \\
{\it{Universit\'e Blaise Pascal,  63177 Aubi\`ere Cedex, France}}  \\
{\tt{E-mail: Fabrice.Blache@math.univ-bpclermont.fr}}              \\
\end{tabular}}
\date{June 2004}

\maketitle

\vskip .5cm

\centerline{\small{\bf Abstract}}

\begin{center}
\begin{minipage}[c]{330pt}
{\small
The problem of finding a martingale on a manifold with a fixed random
terminal value can be solved by considering BSDEs with a generator with
quadratic growth. We study here a generalization of these equations and
we give uniqueness and existence results in two different frameworks,
using differential geometry tools. Applications to PDEs are given,
including a certain class of Dirichlet problems on manifolds.}
\end{minipage}
\end{center}

\bigskip

\section{Introduction}
\label{par1}

\subsection{Martingales and BSDEs on manifolds}
\label{par1.1}
Unless otherwise stated, we shall work on a fixed finite time interval
$[0;T]$; moreover, $(W_t)_{0 \le t \le T}$ will always denote a Brownian
Motion (BM for short) in $\R^{d_W}$, for a positive integer $d_W$.
Moreover, Einstein's summation convention will be used for repeated
indices in lower and upper position.

It is well-known that there is a deep interplay between on the one hand
the probability theory of real martingales and Brownian motion and on the
other hand the theory of PDEs and harmonic functions
$h : M_1 \rightarrow \R$, defined on a manifold $M_1$. For instance,
the Feynman-Kac formula gives a probabilistic interpretation for the
solution of a PDE; besides, the Dirichlet problem for such harmonic
functions $h$ can be solved by considering real martingales with fixed
terminal value.

It is natural to ask whether these links can be generalized to the
nonlinear context of manifolds, i.e. if we replace the vector space $\R$
in the definition of $h$ by a manifold $M$; this is what we now examine.
Suppose that $M$ is a manifold endowed with a connection $\Gamma$; then
one can define the notion of $\Gamma$-martingale on $M$, which generalizes
real local martingales (for an overview of the basic
definitions and properties, see \cite{scm}, \cite{darl82} or \cite{sam}).

In $\R^n$, the problem of finding a martingale $(X_t)_{0 \le t \le T}$ with
terminal value $X_T=U$ consists of solving the Backward Stochastic
Differential Equation (BSDE for short)
$$(E) \left\{
      \begin{array}{l}
          \overline X_{t+dt} = \overline X_t + \overline Z_t d W_t    \\
	  \overline X_T=U,                         \\
      \end{array}
      \right.    $$
where $(\overline Z_t)_{0 \le t \le T}$
is a $\R^{n \times d_W}$-valued progressively measurable process.
With the connection $\Gamma$ on $M$, one can define an exponential mapping
${\rm exp}$ and the equation under infinitesimal form $(E)$ becomes, for
martingales on $M$,
\begin{equation}
(M)_0
\left\{
      \begin{array}{l}
           X_{t+ d t}= {\rm exp}_{X_t}(Z_t d W_t)  \\
           X_T=U                         \\
      \end{array}
      \right.
\end{equation}
where $Z_t \in {\cal L}(\R^{d_W}, T_{X_t}M)$ is now a linear map.
As in the linear context, studying martingales on $M$ (or equivalently
solving BSDE $(M)_0$) allows to solve in a probabilistic way
some nonlinear PDEs and to study harmonic mappings. Let us recall the
definition and some properties of these mappings. \\
A harmonic map $H : M_1 \rightarrow M$ between Riemannian manifolds $M_1$
and $M$ is a smooth map which is a local extremal of the energy
functional
$$\int \Vert {\rm grad} H \Vert^2 d vol$$
where $d vol$ is the Riemannian volume element on $M_1$. \\
A different but equivalent point of view about these mappings is the one
of a system of elliptic PDEs (see \cite{eellslem1}); let us make precise it.
Consider a second-order differential operator ${\cal L}$ without term of
order 0, defined on $M_1$. For $h : M_1 \rightarrow M$, one can define
by means of ${\cal L}$ the tension field of $h$; it is a vector field
along $h$, i.e.
$${\cal L}_M (h) : M_1 \rightarrow TM, \ \ {\cal L}_M(h)(x)
\in T_{h(x)} M.$$
Then the equation ${\cal L}_M (h) = 0$
characterizes ${\cal L}$-harmonic maps (see \cite{eellsamp} and
\cite{eellslem1}, and probabilistic interpretations in the introductions
of \cite{mpl} and \cite{barmart}).
In coordinates $(x^i)$ on $M$ and $(y^\alpha)$ on $M_1$,
this equation can be written as the following system of elliptic PDEs
$$ \forall i, \ \ \Delta_{M_1} \phi^i + \Gamma^i_{jk}(\phi)
    g^{\al \beta}(x) D_\alpha \phi^j D_\beta \phi^k =0$$
with $\Delta_{M_1}$ denoting the Laplace-Beltrami operator and
$(g^{\al \beta})$ the inverse metric tensor on $M_1$. Note that we have
used the summation convention. \\
With the theory of martingales on manifolds, one can solve such a system
of nonlinear elliptic PDEs; for further details, the reader is referred
to \cite{kend90} and \cite{kend94}.

\bigskip
Now we come to the aim of this article, by enlarging the class of
processes studied. Let $(B_t^y)_{0 \le t \le T}$ denote the $\R^d$-valued
diffusion which is the unique strong solution of the following SDE :
\begin{equation}
\left\{
\begin{array}{rcl}
dB_t^y & = & b(B_t^y) dt + \sigma(B_t^y) d W_t  \\
 B_0^y & = & y,
\end{array}
\right.
\label{sde}
\end{equation}
where $\sigma : \R^d \rightarrow \R^{d \times d_W}$ and
$b : \R^d \rightarrow \R^d$ are $C^3$ bounded functions with bounded
partial derivatives of order 1,2 and 3.

On $\R^n$, equation
$(E)$ is a very simple BSDE. A more general form of BSDE on $[0;T]$ is
$$(E+D) \left\{
      \begin{array}{l}
          \overline X_{t+dt} = \overline X_t + \overline Z_t d W_t
            +f(B_t^y, \overline X_t, \overline Z_t) dt  \\
	  \overline  X_T=U,                         \\
      \end{array}
      \right.    $$
studied for instance in \cite{pp90} and \cite{pp92}; such an equation can
be used to solve systems of quasilinear PDEs (see for example \cite{pp92}).
\\
If we try to extend equation $(E+D)$ to manifolds, by combining with
equation $(M)_0$ we get the following equation (under infinitesimal form)
$$
(M+D)_0
\left\{
      \begin{array}{l}
           X_{t+ d t}= {\rm exp}_{X_t}(Z_t d W_t + f(B_t^y, X_t, Z_t) dt)\\
           X_T=U.                         \\
      \end{array}
      \right.
$$
The aim of this work is to study existence and uniqueness of a solution
to the generalized equation $(M+D)_0$. \\
As for martingales, which are linked to harmonic mappings $h$ (i.e. to
the above equation ${\cal L}_M(h) = 0$), this BSDE is related
with a differential operator generalizing the tension field ${\cal L}_M$;
moreover,
in some cases, the class of mappings which solve the new PDE can be
described in terms of the local extrema of another variational problem.
This will be discussed in Section \ref{par5}.
It is known that when $M_1 \subset \R^3$ and $M=S^2$, harmonic mappings
can be used to model the state of equilibrium of liquid crystals (see
the introduction of \cite{helein} for a brief discussion); then mappings
associated to the new variational problem could be used to model the
equilibrium state of a liquid crystal in
an exterior field equal to the drift term $f$ in equation $(M+D)_0$.

\subsection{Setting of the problem}
\label{par1.2}
In the whole paper, we will always suppose that a global system of
coordinates is given on $M$ . Then in these coordinates, we get from
equation $(M)_0$ the following BSDE (see for instance \cite{scm} or the
introduction of \cite{darl95})

$$(M)
\left\{
      \begin{array}{l}
          d X_t = Z_t d W_t - \half \Gamma_{jk}(X_t)
                              ([Z_t]^k \vert [Z_t]^j) d t  \\
	  X_T=U;                                             \\
      \end{array}
    \right.   $$
in this equation, we have used the following notations, which will be
valid throughout the sequel :
$( \cdot \vert \cdot)$ is the usual inner product in an Euclidean space,
the summation convention is used, and $[A]^i$ denotes the $i^{th}$ row
of any matrix $A$; finally,
\begin{equation}
\Gamma_{jk} (x) =
\left(
   \begin{array}{c}
       \Gamma^1_{jk}(x)    \\
       \vdots              \\
       \Gamma^n_{jk}(x)
   \end{array}
\right)
\label{christof}
\end{equation}
is a vector in $\R^n$, whose components are the Christoffel symbols
of the connection. Remark also that here $Z_t$ is a matrix in
$\R^{n \times d_W}$.

In this case,
the classical approach of \cite{pp90} to solve BSDEs with Lipschitz
coefficients fails since there is a quadratic term in $Z_t$ in the
drift (the reader is referred to \cite{pardoux98} or \cite{pengkar97}
for an introduction to the theory of BSDEs in Euclidean spaces).
However, uniqueness and existence results have been obtained
using differential geometry tools, in particular by Arnaudon 
(\cite{arn96}), Darling (\cite{darl95}), Emery(\cite{scm}), Kendall
(\cite{kend90}), Picard (\cite{msc} and \cite{mpl}) or Thalmaier
(\cite{thal96} and \cite{thal96b});
note also the results of Estrade and Pontier (\cite{estradepont}) 
concerning some classes of
Lie groups. Independently of geometric tools, a lot of works have tried
to weaken the Lipschitz assumption : in the one-dimensional case, they
include \cite{kobyl00} (in dimension one, her results are more general
than the ones of this paper because she deals with generators with
quadratic growth), \cite{lepsanm97} or \cite{hamadene96}; in higher
dimensions, we refer the reader for instance to \cite{briandcarm00},
\cite{tang03} (who studies a Ricatti-type BSDE) and \cite{briandhupard03}.
To the best of our knowledge, there is no paper that would include our
results in dimensions greater than one.

\bigskip
Now in our global chart, the equation $(M+D)_0$  becomes
$$(M+D)
\left\{
      \begin{array}{l}
          d X_t = Z_t d W_t + \left( - \half \Gamma_{jk}(X_t)
            ([Z_t]^k \vert [Z_t]^j) + f(B_t^y,X_t,Z_t) \right) d t  \\
	  X_T=U                                               \\
      \end{array}
    \right.   $$
(the same notation will be used to denote the $TM$-valued function $f$
and its image in local coordinates).
The process $X$ will take its values in a compact set, and a solution of
equation $(M+D)_0$ will be a pair of processes
$(X,Z)$ in $M \times (\R^{d_W} \otimes TM)$ such that $X$ is continuous and
$\esp \left( \int_0^T \Vert Z_t \Vert_r^2 dt \right) < \infty$
($\Vert \cdot \Vert_r^2$ is a Riemannian norm; see below).
If we consider a global system of coordinates on an open set $O$ of
$\R^n$, it corresponds to processes $(X,Z)$ in $O \times \R^{n d_W}$,
such that $X$ is in a compact set and
$\esp \left( \int_0^T \Vert Z_t \Vert ^2 dt \right) < \infty$, solving
equation $(M+D)$.

Two different cases will be considered here : firstly when the drift
$f$ does not depend on $z$, and secondly the case of a general $f$
in nonpositive curvatures. In the two cases, a Riemannian structure is
fixed on $M$; in the former case, the connection may be independent of
this Riemannian structure, while in the latter, only the Levi-Civita
connection associated will be used. Note that the case for a general
$f$ with $K>0$ involves more technical calculations; it will appear
elsewhere.

We first give in Section \ref{par2} mild generalizations of well-known
results, concerning the geometry of the manifold and a characterization
of the solutions in the $z$-independent case by means of convex
functions.
In Section
\ref{par3}, we study the uniqueness problem. It is solved by generalizing
to our context two methods : on the one hand, Emery's idea, used in
\cite{scm} and \cite{kend90}; on the other hand, the work of Picard
(\cite{mpl}). We obtain Theorems \ref{unicit1} and \ref{unicit2}.

Section \ref{par4} is devoted to proving the existence of a solution of
equation $(M+D)$. The main arguments are to exhibit a solution for
"simple" terminal values (based on a strong bound on the process $(Z_t)$
in Subsection \ref{par4.3}) and to
solve the equation for any terminal value using approximation procedures
(Subsections \ref{par4.1} and \ref{par4.6}).
We need for the proof an additional (and necessary in fact) condition
on the drift $f$~: it
is supposed to point outward on the boundary of the set on which we work.
We give in Subsection \ref{par1.4} the main result (Theorem
\ref{existence2}) which sums up the results obtained.
In Section \ref{par5}, we extend the results to random time intervals
$[0; \tau]$, where $\tau$ is successively a bounded stopping time (Theorem
\ref{existence3}) and a stopping time verifying an exponential
integrability condition 
(Theorem \ref{existence4}); then to conclude this paper, we give some
generations and applications to the theory of PDEs, as well as
the variational problem related to equation $(M+D)_0$.

The uniqueness part, as well as the applications to PDEs, are mainly
adaptations of procedures already used; on the contrary, the approach
for the existence seems to be novel.

\subsection{Notations and hypothesis}
\label{par1.3}
In all the article, we suppose that a filtered probability space
$(\Omega, {\cal F}, P, ({\cal F}_t)_{0 \le t \le T})$ is given on which
$(W_t)_t$ denotes a $d_W$-dimensional BM. Moreover, we
always deal with a complete Riemannian manifold  $M$ of dimension $n$,
endowed with a linear symmetric (i.e. torsion-free) connection whose
Christoffel symbols $\Gamma^i_{jk}$ are smooth; the connection does not
depend {\it a priori} on the Riemannian structure.

On $M$, $\delta$ denotes the Riemannian distance;
$\vert u \vert_r$ is the Riemannian norm for a tangent vector $u$
and $\vert u' \vert$ the Euclidean norm for a vector $u'$ in $\R^n$.
If $h$ is a smooth real function defined on $M$ and $u_1, u_2$ are tangent
vectors at $x$, the differential of $h$ is denoted by $D h(x)<u>$ or
$h'(x)<u>$;
the Hessian $\Hess h(x)$ is a bilinear form the value of which is denoted
by $\Hess h(x)<u_1,u_2>$.

For $\beta \in \N^*$, we say that a function is $C^\beta$ on a closed set
$F$ if it is $C^\beta$ on an open set containing $F$.
For a matrix $z$ with $n$ rows and $k$ columns, $^t z$ denotes
its transpose,
$$\Vert z \Vert
= \sqrt{ \hbox{Tr} (z {}^t z)}
= \sqrt{\sum_{i=1}^k \vert [{}^t z]^i \vert ^2}$$
(${\rm Tr}$ is the trace of a square matrix) and
$\Vert z \Vert_r = \sqrt{\sum_{i=1}^k \vert [{}^t z]^i \vert_r ^2}$
where the columns of $z$ are considered as tangent vectors.
The notation $\Psi(x,x') \approx \delta(x,x')^\nu$ means that there is
a constant $c>0$ such that
$$ \forall x,x', \
\frac{1}{c} \ \delta(x,x')^\nu \le \Psi(x,x') \le c \ \delta(x,x')^\nu.$$
Finally, recall that a real function $\chi$ defined on $M$ is said to be
convex if for any $M$-valued geodesic $\gamma$, $\chi \circ \gamma$ is
convex in the usual sense (if $\chi$ is smooth, this is equivalent to
require that $\Hess \chi$ be nonnegative).

\medskip
Before the general framework, let us give some additional notations
which are specific to the Levi-Civita connection.
In this case, we always assume that the injectivity radius $R$ of $M$
is positive and that its sectional curvatures are bounded above; we let
$K$ be the smallest nonnegative number dominating all the sectional
curvatures.
For the distance function $\delta$, if $\tilde x =(x,x')$ is a point
and $u, \overline u$ are tangent vectors at $x'$, we consider the
partial derivatives (when they are defined)
$$ \delta'_2(\tix)<u> = \delta'(\tix)<(0,u)>$$
and
$$ {\rm Hess}_{22} \  \delta (\tix)<u,\overline u> =
     \Hess \delta (\tix)<(0,u),(0,\overline u)>.$$

\noi Now let us recall from \cite{kend90}
the definition of a regular geodesic ball. A closed geodesic ball $\cal B$
of radius $\rho$ and center $p$ is said to be regular if

\noi (i) $\rho \sqrt K < \frac{\pi}{2}$

\noi (ii) the cut locus of $p$ does not meet $\cal B$.
\\
For an introductory course in Riemannian
geometry, the reader is referred to \cite{booth} and for further facts
about curvature, to \cite{lee}.

\bigskip
Throughout this article, we consider an open set $\omega \not= \emptyset$
relatively compact in an open subset $O$ of $M$, such that there is a
unique $O$-valued geodesic between any two points of $O$; $O$ is also
supposed to be relatively compact in a local chart, so that it provides a
system of coordinates (global on $\omb$); it will be as well
considered as a subset of $\R^n$. We suppose that there exists a
nonnegative, smooth and convex function $\Psi$ on the product
$\omb \times \omb$ (i.e. convex on an open set containing this set) which
vanishes only on the diagonal $\Delta=\{(x,x)/x \in \omb \}$ ($\omb$ is
said to have $\Gamma$-convex geometry); besides, we suppose that
$\Psi \approx \delta^p$ for a $p \ge 2$ (note that since $\Psi$ is
smooth, $p$ is an even integer).
In fact, we take for $\omb$ a sublevel set of a smooth convex function
$\chi$ defined on  $O$ : $\{ \chi \le c \}$ (note that this hypothesis and
the existence of $\Psi$ guarantee the existence and uniqueness of a
$\omb$-valued geodesic between any two points of $\omb$).
\\
Emery has shown (see Lemma (4.59) of \cite{scm}) that, in the case of a
general connection $\Gamma$, any point of $M$ possesses a neighbourhood
with $\Gamma$-convex geometry; when the Levi-Civita connection is used,
this is true for a regular geodesic ball (see \cite{kend91}).

Finally we always assume two hypothesis on $f$ :
$$  \exists L>0 , \ \forall b,b' \in \R^d,
  \forall (x,z) \in O \times {\cal L}(\R^{d_W}, T_xM),
  \forall (x',z') \in O \times {\cal L}(\R^{d_W}, T_{x'}M),
$$
\begin{eqnarray}
\left\vert \overset{x'}{\underset{x}{\Vert}} f(b,x,z) - f(b',x',z')
\right\vert_r \le L \Bigg( ( \vert b-b' \vert
& + &\delta(x,x')) (1 + \Vert z \Vert_r + \Vert z' \Vert_r)  \nonumber  \\
& + & \left\Vert \overset{x'}{\underset{x}{\Vert}} z - z' \right\Vert_r
                 \Bigg)
\label{lip}
\end{eqnarray}
and
\begin{equation}
\exists L_2>0 , \exists x_0 \in O, \forall b \in \R^d,
\vert f(b,x_0,0) \vert_r \le L_2.
\label{upperboundf}
\end{equation}
The first one is a "geometrical" Lipschitz condition on $f$. This special
form is needed to get an expression which is invariant under changes of
coordinates. We will see that later, in \eqref{lip2}. Remark that, in the
$z$-independent case, it just means that $f$ is Lipschitz with
respect to the first two variables; otherwise,
$$\overset{x'}{\underset{x}{\Vert}}  z$$
 denotes the Riemannian parallel
transport along the unique geodesic between $x$ and $x'$.
The second one means that $f$ is bounded with respect to the first
argument.
Remark that these conditions also imply the boundedness of $f$ if it does
not depend on the variable $z$.

Note to end this part that the same letter $C$ will often stand for
different constant numbers.

\subsection{The main result}
\label{par1.4}
Before achieving calculations, we give the main theorem of the article.
Let us first introduce a technical but natural hypothesis, which we will
make explicit in Subsection \ref{par4.6}~:
$$(H) \ \ f \hbox{ is pointing outward on the boundary of } \omb.$$
Then we can state :
\begin{stheorem}
\label{existence2}
We consider the BSDE $(M+D)$
with terminal random variable $U \in \omb= \{ \chi \le c \}$, where $\omb$
satisfies the above conditions. If $f$
verifies conditions \eqref{lip}, \eqref{upperboundf} and $(H)$, and
if $\chi$ is strictly convex (i.e. $\Hess \chi$ is positive definite),
then
\begin{item}
(i) If f does not depend on $z$, the BSDE has a unique solution
$(X_t,Z_t)_{0 \le t \le T}$ such that $X$ remains in $\omb$.
\end{item}
\begin{item}
(ii) If $M$ is a Cartan-Hadamard manifold and the Levi-Civita connection
is used, then the BSDE has yet a unique solution
$(X_t,Z_t)_{0 \le t \le T}$ with $X$ in $\omb$.
\end{item}
\end{stheorem}
In particular, if the Levi-Civita connection is used, a "good" example of
domain on which existence and uniqueness hold is a regular geodesic ball.

In Section \ref{par5}, we will extend this theorem to random time
intervals $[0;\tau]$ (instead of $[0;T]$), for stopping times $\tau$ which
are bounded, or verify the exponential integrability condition~:
\begin{equation}
  \exists \rho > 0 :  \esp (e^{\rho \tau})<\infty.
\label{integst}
\end{equation}
In the former case, Theorem
\ref{existence2} goes the same, while in the latter the constants $L$ and
$L_2$ in \eqref{lip} and \eqref{upperboundf} are furthermore required to
be small with respect to the constant $\rho$ in \eqref{integst}.

\bigskip
{\it Acknowledgements : } The author would like to thank his supervisor
Jean Picard for his help and his relevant advice, and the referees for
their suggestions to improve a first version.

\section{Preliminary results}
\label{par2}

We first recall elementary results about It\^o's formula and parallel
transport.
Then we give some geometrical
estimates for the distance function on $M \times M$ and characterize
solutions of the equation $(M+D)$ using convex functions, but only when
the drift $f$ does not depend on $z$. As underlined in the introduction,
these results are just mild generalizations of well-known results of
\cite{scm} and \cite{mpl}.

In this section, the covariant derivative of a vector field $z_t$ along
a curve $\gamma_t$ will be denoted $\nabla_{\dot \gamma_t} z_t$.

\subsection{It\^o's formula on manifolds}

Consider two solutions $(X^1,Z^1)$ and $(X^2,Z^2)$ of equation
$(M+D)$ with terminal values $U^1$ and $U^2$, such that $X^1$
and $X^2$ remain in $O$. Let
$$\tilde X =(X^1,X^2) \ \ \hbox{ and } \ \
    \tilde Z = \left(
                 \begin{array}{c}
                    Z^1   \\
	            Z^2
                  \end{array}
                \right);$$
then It\^o's formula with the function $\Psi$ is written
	
\begin{eqnarray}
\Psi (\tilde X_t) - \Psi (\tilde X_0)
    &=&  \int_0^t D \Psi (\tilde X_s)
           \left( \tilde Z_s d W_s \right)             \nonumber   \\
    & & +  \int_0^t D \Psi (\tilde X_s)
            \left(
                 \begin{array}{c}
          f(B_s^y,X^1_s,Z^1_s)- \half \Gamma_{jk}(X^1_s)
                                 ([Z^1_s]^k \vert [Z^1_s]^j)  \\
          f(B_s^y,X^2_s,Z^2_s)- \half \Gamma_{jk}(X^2_s)
                                 ([Z^2_s]^{k} \vert [Z^2_s]^{j})
                 \end{array}
            \right)  ds                                \nonumber   \\
    & & +  \half \int_0^t  {\rm Tr}  \left(
                   {}^t \tilde Z_s D^2 \Psi(\tilde X_s) \tilde Z_s
                                          \right) ds    \nonumber  \\
    &=& \int_0^t D \Psi (\tilde X_s)
            \left( \tilde Z_s d W_s \right)            \nonumber   \\
    & & + \half \int_0^t  \left(
      \sum_{i=1}^{d_W} {}^t [{}^t \tilde Z_s]^i \Hess \Psi (\tilde X_s)
                  [{}^t \tilde Z_s]^i \right) ds         \nonumber \\
    & & + \int_0^t  D \Psi(\tilde X_s)
                                   \left(
                                 \begin{array}{c}
                                     f(B_s^y,X^1_s,Z^1_s)  \\
                                     f(B_s^y,X^2_s,Z^2_s)
                                 \end{array}
                                   \right) ds
\label{ito1}
\end{eqnarray}
(remember notation \eqref{christof} and that $[{}^t \tilde Z_s]^i$
denotes the $i^{th}$ column of the
matrix $\tilde Z$; it is a vector in $\R^{2n}$).
Moreover, for a smooth function $h$ on $O$ and a solution $(X,Z)$ of
$(M+D)$, we get a similar formula, replacing $\tilde X$ by $X$ and
$\tilde Z$ by $Z$.

\subsection{Two inequalities}

We first give an equivalence result between the Euclidean and Riemannian
norms; it follows easily from the relative compactness of $O$, considered
as a subset of $\R^n$ (in particular, this means that we can identify
each tangent space with $\R^n$).

\begin{slemma}
\label{equival}
There is a $c>0$ such that for any
$(x,z) \in O \times T_xM = O \times \R^n$,
$$ \frac{1}{c} \vert z \vert \le \vert z \vert_r \le c \vert z \vert.$$
\end{slemma}
Lemma \ref{equival} will often be useful in the sequel.

\begin{sprop}
\label{proptp}
The Levi-Civita connection is used. There is a $C>0$ such that for every
$(x,x') \in O \times O$
and  $(z,z') \in T_xM \times T_{x'}M$, we have

\begin{equation}
\left\vert \overset{x'}{\underset{x}{\Vert}} z - z' \right\vert_r
 \le C
  \left( \vert z-z' \vert + \delta(x,x')(\vert z \vert + \vert z' \vert)
       \right)
 \label{tp1}
\end{equation}
and
\begin{equation}
  \vert z-z' \vert \le C
\left( \left\vert \overset{x'}{\underset{x}{\Vert}}
                                    z - z' \right\vert_r
  + \delta(x,x')(\vert z \vert_r + \vert z' \vert_r)
       \right).
 \label{tp2}
\end{equation}

\end{sprop}

{\bf Remark 1.} In fact, by Lemma \ref{equival}, we can use any of the
two norms (except for $\vert z-z' \vert$ which is necessarily the
Euclidean norm).

\begin{demo}
It is sufficient to prove
\begin{equation}
\forall (x,z), (x',z') \in O \times \R^{n},
\left\vert \overset{x'}{\underset{x}{\Vert}} z - z \right\vert
 \le C \delta(x,x') \vert z \vert.
\label{tp3}
\end{equation}
Indeed, this implies

$$\left\vert \overset{x'}{\underset{x}{\Vert}} z - z' \right\vert_r
  \le C_1 \left( \vert z-z' \vert +
\left\vert \overset{x'}{\underset{x}{\Vert}} z - z \right\vert \right)
  \le  C ( \vert z-z' \vert + \delta(x,x') \vert z \vert)$$

and

$$ \vert z-z' \vert \le
   \left\vert \overset{x'}{\underset{x}{\Vert}} z - z \right\vert
 + \left\vert \overset{x'}{\underset{x}{\Vert}} z - z' \right\vert
 \le C \left( \delta(x,x') \vert z \vert +
  \left\vert \overset{x'}{\underset{x}{\Vert}} z-z'\right\vert \right).$$

So let us prove \eqref{tp3}~:
let $\gamma$ be the geodesic such that $\gamma(0)=x$ and $\gamma(1)=x'$,
and $z(t)$ the parallel transport of $z$ along $\gamma$ :
$$\forall t \in [0;1], \ \
   z(t)=\overset{\gamma_t}{\underset{x}{\Vert}} z.$$
In local coordinates, the equation $\nabla_{\dot \gamma_t} z(t)=0$ gives
for every $k$

$$ \dot z^k(t) + \Gamma_{jl}^k(\gamma_t) \dot \gamma_t^j z^l(t)=0.$$
Moreover, $\vert \dot \gamma_t \vert_r=\delta(x,x')$ so
$\vert \dot \gamma_t \vert \le \tilde C \delta(x,x')$ and
\begin{eqnarray*}
\left\vert \overset{x'}{\underset{x}{\Vert}} z - z \right\vert ^2
& =  & \vert z(1)-z(0) \vert ^2                     \\
& =  & \sum_k \vert z^k(1)-z^k(0) \vert ^2          \\
& =  & \sum_k \left\vert \int_0^1 \dot z^k(t) dt \right\vert ^2 \\
& =  & \sum_k \left\vert \int_0^1 \Gamma_{jl}^k(\gamma_t) \dot \gamma_t^j
           z^l(t) dt \right\vert ^2                       \\
&\le & C_1 \delta(x,x')^2 \sum_l \int_0^1 \vert z^l(t) \vert ^2 dt  \\
&\le & C \delta(x,x')^2 \vert z \vert ^2.
\end{eqnarray*}
The last inequality comes from the equivalence on $O$ of the Riemannian
and Euclidean norms, and the fact that $\vert z(t) \vert_r$ is
constant by definition of $z(t)$.
The proof is completed.
\end{demo}

As a consequence, on the relatively compact set $O \subset \R^n$,
\eqref{lip} becomes

$$  \exists L'>0 , \ \forall b,b' \in \R^d, \
  \forall (x,z), (x',z') \in O \times \R^{n d_W},
$$
\begin{equation}
\vert  f(b,x,z) - f(b',x',z') \vert \le L'
( (\vert b-b' \vert  +  \vert x - x' \vert)
        (1 + \Vert z \Vert + \Vert z' \Vert) +  \Vert z - z' \Vert ).
\label{lip2}
\end{equation}

{\bf Remark 2.} In particular, \eqref{lip2} is verified for a drift $f$
that is Lipschitz in $(b,x,z)$ in $O$ (the Lipschitz property of $f$ is
not necessarily preserved by a change of coordinates, but \eqref{lip2}
is).

{\bf Remark 3.} If $f$ does not depend on $z$, it just means that $f$ is
Lipschitz in $(b,x)$.

\subsection{Estimates of the derivatives of the distance}
\label{estimdist}

This paragraph is based on Section 1 of \cite{mpl}. The
connection used is Levi-Civita's one.

The geodesic distance $(x,x') \mapsto \delta(x,x')$
is defined on $M \times M$ and is smooth except on the cut locus and the
diagonal $\{x=x'\}$. We want
to estimate its first and second derivatives when $M \times M$ is endowed
with the product Riemannian metric. If $\tix=(x,x')$ is a point which is
not in the cut locus or the diagonal, there exists a unique minimizing
geodesic $\gamma(t)$, $0 \le t \le 1$, from $x$ to $x'$. If $u_t$ is a
vector of $T_{\gamma(t)}M$, we
can decompose $u_t$ as $v_t+w_t$, where $v_t$ is the orthogonal projection
of $u_t$ on $\dot \gamma(t)$; the vectors $v_t$ and $w_t$ are respectively
called
the tangential and orthogonal components of $u_t$. If $u=(u_0,u_1)$ is a
vector of $T_{\tix} (M \times M)$, $(v_0,v_1)$ and $(w_0,w_1)$ are also
called its tangential and orthogonal components.

\begin{slemma}
\label{derpsi}
Let $\tix$ be a point of $M \times M$ which is not in the cut locus
or the diagonal.
Let $u$ be a vector of $T_{\tix}(M \times M)$ and let $v$ and $w$ be its
tangential and orthogonal components. Then
\begin{equation}
  \vert \delta'(\tix)<u> \vert =
  \left\vert \overset{x'}{\underset{x}{\Vert}} v_0 - v_1 \right\vert _r ;
\label{der1}
\end{equation}
if moreover $K=0$ (i.e. the sectional curvatures are nonpositive), then
\begin{equation}
 \Hess \delta(\tix)<u,u> \ge
  \frac{1}{\delta(\tix)} \left\vert \overset{x'}{\underset{x}{\Vert}}
        w_0 - w_1  \right\vert_r^2.
\label{der2}
\end{equation}
\end{slemma}

\begin{demo}
Let $J_v(t)$ (resp. $J_w(t)$) be the tangential (resp. normal) Jacobi
field along $\gamma(t)$ satisfying $J_v(0)=v_0$ and $J_v(1)=v_1$ (resp.
$J_w(0)=w_0$ and $J_w(1)=w_1$). From (1.1.5) and
(1.1.7) of \cite{mpl}, we have
\begin{equation}
\delta'(\tix)<u> =
   \frac{(\dot \gamma(t) \vert \nabla_{\dot \gamma(t)} J_v(t))}
        {\vert \dot\gamma(t) \vert _r}
\label{der1b}
\end{equation}
and
\begin{equation}
 \Hess \delta(\tix)<u,u> \ge \frac{1}{\delta(\tix)} \int_0^1
       \vert \nabla_{\dot \gamma(t)} J_w(t) \vert_r^2 dt
    - K \delta(\tix) \int_0^1 \vert J_w(t) \vert_r^2 dt.
\label{der2b}
\end{equation}
It is easy to see with the Jacobi equation that
$J_v(t)=(At+B) \dot \gamma(t)$. Then
$\nabla_{\dot \gamma(t)} J_v(t) = A \dot \gamma(t)$ and
the limit conditions $J_v(0)=v_0= \alpha \dot \gamma(0)$ and
$J_v(1)=v_1= \beta \dot \gamma(1)$ imply $A= \beta -\alpha$. Hence
\eqref{der1b} gives
$$\vert \delta'(\tix)<u> \vert =
              \vert A \vert \cdot \vert \dot \gamma(t)\vert_r =
          \vert \beta -\alpha \vert \delta(\tilde x) =
         \left\vert \overset{x'}{\underset{x}{\Vert}} v_0 - v_1
     \right\vert_r.$$
This is \eqref{der1}.

Moreover, let us write $J_w(t)=\sum_i u_i(t) E_i(t)$ where
$\{ E_i(t) \}_i$ is a parallel orthonormal frame along $\gamma$. Its
covariant derivative is
$\nabla_{\dot \gamma(t)} J_w(t)= \sum_i \dot u_i(t) E_i(t)$ and
$ \overset{x'}{\underset{x}{\Vert}} J_w(0)= \sum_i u_i(0) E_i(1)$. Then
\begin{eqnarray*}
\int_0^1 \vert \nabla_{\dot \gamma(t)} J_w(t) \vert_r^2 dt
& = & \int_0^1 \sum_i \vert \dot u_i(t) \vert ^2 dt  \\
&\ge& \sum_i \left( \int_0^1 \dot u_i(t) dt \right) ^2
  = \left\vert \overset{x'}{\underset{x}{\Vert}} J_w(0) - J_w(1)
       \right\vert_r ^2.
\end{eqnarray*}
Now this inequality together with \eqref{der2b} and the nonpositivity
of the sectional curvatures give \eqref{der2}.
\end{demo}

Then we have the following estimate~:
\begin{sprop}
If $K=0$ and $\tix$ is not in the cut locus, then
\begin{equation}
           \Hess (\half \delta^2)(\tix)<u,u> \ge
            \left\vert \overset{x'}{\underset{x}{\Vert}}
                u_0 - u_1 \right\vert_r ^2.
\label{min}
\end{equation}
\end{sprop}

\begin{demo}
If $\tilde x=(x,x')$ is not on the diagonal, we recall the classical
formula
$$\Hess (\half \delta^2)(\tix)<u,u>=\delta(\tix) \cdot
\Hess \delta(\tix)<u,u> + ( \delta'(\tix)<u>)^2.$$
This formula and estimates \eqref{der1} and \eqref{der2} imply the
proposition for $x \not= x'$ since the two terms
$\overset{x'}{\underset{x}{\Vert}} v_0 - v_1$ and
$\overset{x'}{\underset{x}{\Vert}} w_0 - w_1$
are orthogonal for
the Riemannian norm. The case $x=x'$ follows by continuity since
$\delta^2$ is smooth on a neighbourhood of the diagonal.
\end{demo}

\subsection{A characterization of the solutions of equation $(M+D)$ when
the drift  $f$ does not depend on $z$}

We give here a generalization of a well-known result (see
(4.41)(ii) in \cite{scm}) which roughly says that a continuous
$M$-valued process $(Y_t)$ is a $\Gamma$-martingale if and only if its
image under convex functions is a real local submartingale. In this
paragraph, the filtration used is the natural one of $(W_t)_t$.

\begin{sprop}
\label{charact}
Suppose that the drift $f$ does not depend on $z$. Then every point
$p$ of $M$ has an open neighbourhood $O_p$, included in a local chart,
with the following property :

A pair of processes $(X,Z)$ (with $X$ continuous, adapted and
$O_p$-valued) is a solution of $(M+D)$ iff for every convex function
$\xi: O_p \rightarrow \R$,
$\xi(X_t) - \int_0^t D \xi (X_s) \cdot f(B_s^y,X_s) ds$
is a local submartingale.
\end{sprop}

The proof given here is just an adaptation of Emery's one; first we recall
Lemma (4.40) of \cite{scm}.

\begin{slemma}
\label{approxaff}
On $M$, let $\xi$ be a smooth function. Every point of $M$ has an open
neighbourhood $O_p$ depending on $\xi$ with the following property~:
For every $\ep>0$ and $a \in O_p$, there is a convex function
$h^a_\ep : O_p \rightarrow \R$ such that $(a,x) \mapsto h_\ep^a(x)$ is
smooth on $O_p \times O_p$, $h_\ep^a(a)=0$, $Dh_\ep^a(a)=D\xi(a)$ and
$\Hess h^a_\ep(a)=\ep g(a)$, where $g$ represents the metric.
\end{slemma}

Now we complete the proof of Proposition \eqref{charact}~:

\noindent
The "only if" part is just a consequence of It\^o's formula (similar to
\eqref{ito1})
applied to $\xi(X_t)$~: as $\Hess \xi$ is nonnegative by convexity,
$\xi(X_t) - \int_0^t D \xi (X_s) \cdot f(B_s^y,X_s) ds$ is indeed a local
submartingale.

For the "if" part, notice first that around $p$ there is a system $(x^i)$
of local coordinates
that are convex functions (if $(y^i)$ are any local coordinates with
$y^i(p)=0$, then take $x^i=y^i+ c \sum_j (y^j)^2$ for $c>0$ large enough).
Choose $O_p$ relatively compact in the domain of such a local chart and in
an open set on which Lemma \eqref{approxaff} holds for $\xi=x^i$ and
$\xi=-x^i$ (denote by $S$ this set of $2n$ functions).

For a continuous adapted $O_p$-valued process $X$, suppose that
$h \circ X - \int D h (X) \cdot f(B^y,X) dt$ is a local submartingale for
every convex $h$ on $O_p$. Taking
first for $h$ the global (on $O_p$) coordinates $(x^i)$ shows that each
$x^i \circ X$ is a real semimartingale, so $X$ is a semimartingale. For
every fixed $\xi \in S$, it is sufficient to prove that
$ \xi \circ X - \half \int (\Hess \xi)_{ij} d <X^i,X^j>
   - \int D_i \xi (X) \cdot f^i(B^y,X) dt$
is a local submartingale; for then replacing $\xi$ by $- \xi$ shows that
(remember that $D_{ij} x^k=0$ and $D_i x^k=1$ if $i=k$ and $0$ otherwise)
$$x^k \circ X + \half \int \Gamma^k_{ij}(X) d <X^i,X^j>
    - \int f^k(B^y,X) dt$$
is a local martingale for each $k$. But the theorem of representation of
local martingales in Brownian filtrations allows to write this local
martingale explicitly as $\int Z_t d W_t$; thus
$d <X^i,X^j>_t = ([Z_t]^i \vert [Z_t]^j) dt$ and $(X,Z)$ solves equation
$(M+D)$.

By the choice of $O_p$, given any $\ep>0$ we are provided with functions
$h_\ep^a$ associated to $\xi$ as in Lemma \eqref{approxaff}. Call $\sigma$
the $p^{th}$ dyadic subdivision of the time axis
($\sigma= \{ t_k=\frac{k}{2^p} \ : \ k,p \in \N \}$) and let for
$t \in [t_k; t_{k+1}[$, $\rho(t)=t_k$ and
\begin{eqnarray*}
 S_t^\sigma = \sum_{l<k} \bigg(  h_\ep^{X_{t_l}}(X_{t_{l+1}})
& - &\int_{t_l}^{t_{l+1}} D h_\ep^{X_{t_l}}(X_u) \cdot f(B_u^y,X_u)
                                                            du \bigg) \\
& + & h_\ep^{X_{t_k}}(X_t)
   - \int_{t_k}^t D h_\ep^{X_{t_k}}(X_u) \cdot f(B_u^y,X_u) du.
\end{eqnarray*}
As each $h_\ep^a$ is convex,
$h_\ep^{X_{t_k}}(X_t) - \int_{t_k}^t D h_\ep^{X_{t_k}}(X_u)
\cdot f(B_u^y,X_u) du$
is a submartingale in the interval $[t_k;t_{k+1}]$, and $S^\sigma$ is a
continuous submartingale. Using the coordinates $(x^i)$, write
$$ d S_t^\sigma = D_i h_\ep^{X_{t_k}}(X_t) d X^i_t
     + \half D_{ij} h_\ep^{X_{t_k}}(X_t) d<X^i,X^j>_t
     - D_i h_\ep^{X_{t_k}}(X_t) f^i(B_t^y,X_t) dt.$$
So, if $X^i$ is decomposed into $N^i+A^i$ (i.e. local martingale $+$
bounded variation part),
$$ d S_t^\sigma
  - \left(
     D_i h_\ep^{X_{t_k}}(X_t) d A^i_t
     + \half D_{ij} h_\ep^{X_{t_k}}(X_t) d<X^i,X^j>_t
     - D_i h_\ep^{X_{t_k}}(X_t) f^i(B_t^y,X_t) dt
    \right)$$

is a local martingale; hence the process

\begin{eqnarray*}
B_t^\sigma = \int_0^t D_i h_\ep^{X_{\rho(s)}}(X_s) dA^i_s
& + & \half \int_0^t D_{ij} h_\ep^{X_{\rho(s)}}(X_s) d<X^i,X^j>_s \\
& - & \int_0^t D_i h_\ep^{X_{\rho(s)}}(X_s) \cdot f^i(B_s^y,X_s) ds
\end{eqnarray*}

is increasing. The estimates (using relative compactness of $O_p$)

\begin{eqnarray*}
\vert D_i h _\ep^u(v) -  D_i h _\ep^v(v) \vert
& \le & C \vert u-v \vert                                \\
\vert D_{ij} h _\ep^u(v) -  D_{ij} h _\ep^v(v) \vert
& \le & C \vert u-v \vert
\end{eqnarray*}

 and the convergence of $X_{\rho(t)}$ to $X_t$
when $p$ goes to infinity (i.e. when $\sigma$ becomes finer), yield a
dominated convergence of $D_i h_\ep^{X_{\rho(s)}}(X_s)$ to
$D_i h_\ep^{X_s}(X_s)=D_i \xi (X_s)$ (since $dh_\ep^a(a)=d \xi(a)$) and
of $D_{ij} h_\ep^{X_{\rho(s)}}(X_s)$ to
$\Gamma^k_{ij}(X_s) D_k \xi(X_s) + \ep g_{ij}(X_s)$ (since
$\Hess h^a_\ep(a)=\ep g(a)$). Hence $B^\sigma$ has a limit, equal to

\begin{eqnarray*}
 \int D_i \xi(X) d A^i
& + &\half \int \Gamma^k_{ij}(X) D_k \xi (X) d<X^i,X^j>     \\
& + &\half \ep \int g_{ij}(X) d<X^i,X^j>
       - \int D_i \xi (X) \cdot f^i(B^y,X) dt,
\end{eqnarray*}

that is an increasing process too. Letting now $\ep$ tend to zero,
$$ J =\int D_i \xi(X) d A^i + \half \int \Gamma^k_{ij}(X) D_k \xi (X)
   d<X^i,X^j> - \int D_i \xi (X) \cdot f^i(B^y,X) dt$$

is also increasing, and (remember that $D_{ij} \xi =0$)
\begin{eqnarray*}
 \xi \circ X - \xi \circ X_0 - \half \int (\Hess \xi)_{ij}
        d <X^i,X^j>
&   &                                                               \\
- \int D_i \xi (X) \cdot f^i(B^y,X) dt
& = & \int D_i \xi(X) d N^i + J
\end{eqnarray*}
is a local submartingale, as was to be proved.
\qed \bigskip


\section{The uniqueness property}
\label{par3}
In the first paragraph, we set the problem and exhibit the sum \eqref{pos},
whose nonnegativity suffices to have uniqueness.
Then, we give a useful estimate and derive the result in the two cases
considered in this paper. The calculus is rather longer if the drift $f$
depends on $z$, for we have to prove exponential integrability.

\subsection{The general method}
Consider two solutions $(X_t,Z_t)_t$ and $(X'_t,Z'_t)_t$ of $(M+D)$ such
that $X$ and $X'$ remain in $\omb$ and $X_T=Y_T=U$ (we will often write
"$\omb$-valued solutions of $(M+D)$").
Let
$$\tilde X_s = (X_s,X'_s) \ \ \hbox{ and }
   \tilde Z_s=  \left(
          \begin{array}{c}
             Z_s  \\
             Z'_s
          \end{array}  \right).$$
In the martingale case ($f=0$), It\^o's formula \eqref{ito1} and the
convexity of $\Psi$ ensure that
the process $(\Psi (\tilde X_t))_t$ is a submartingale. The other
properties of $\Psi$ (see the introduction) then imply that $\tilde X$
remains in the diagonal $\Delta$, therefore the uniqueness required.
This is Emery's method (see Corollary (4.61) in \cite{scm}).

For our purpose (i.e. $f$ does not vanish identically), we want to keep
the submartingale property and therefore control the integral involving
$f$ in \eqref{ito1}. The idea is to study, rather than
$(\Psi (\tilde X_t))_t$, the new process
$(\exp (A_t) \Psi (\tilde X_t))_t$
where 
$$
A_t= \lambda t + \mu \int_0^t (\Vert Z_s \Vert_r  +
                           \Vert Z'_s \Vert_r ) ds,
$$
for appropriate nonnegative constants $\lambda$ and $\mu$. Apply It\^o's 
formula to obtain

\begin{eqnarray}
e^{A_t}  \Psi (\tilde X_t) - \Psi (\tilde X_0)
& = &\int_0^t e^{A_s} d (\Psi (\tilde X_s))
+ \int_0^t e^{A_s} (\lambda + \mu (\Vert Z_s \Vert_r + \Vert Z'_s \Vert_r))
                       \Psi (\tilde X_s) ds               \nonumber \\
& = &\int_0^t e^{A_s} D \Psi (\tilde X_s)
          \left( \tilde Z_s d W_s \right)                 \nonumber \\
&   & + \half \int_0^t e^{A_s}
          \left(  \sum_{i=1}^{d_W} {}^t [{}^t \tilde Z_s]^i
   \Hess \Psi (\tilde X_s) [{}^t \tilde Z_s]^i \right) ds \nonumber \\
&  & + \int_0^t e^{A_s} D \Psi(\tilde X_s)
                                   \left(
                                 \begin{array}{c}
                                     f(B_s^y,X_s,Z_s)  \\
                                     f(B_s^y,X'_s,Z'_s)
                                  \end{array}
                                   \right) ds             \nonumber \\
&  & + \int_0^t e^{A_s} \Psi (\tilde X_s) (\lambda + \mu (\Vert Z_s \Vert_r
     + \Vert Z'_s \Vert_r)) ds.
\label{ito3}
\end{eqnarray}

It is clear that the submartingale property will be preserved if we
show the nonnegativity of the sum
\begin{eqnarray}
\half \sum_{i=1}^{d_W} {}^t [{}^t \tilde Z_t]^i
        \Hess \Psi (\tilde X_t) [{}^t \tilde Z_t]^i
& + &  D \Psi(\tilde X_t)
                                   \left(
                                 \begin{array}{c}
                                     f(B_t^y,X_t,Z_t)  \\
                                     f(B_t^y,X'_t,Z'_t)
                                  \end{array}
                                   \right)             \nonumber     \\
& + & (\lambda + \mu (\Vert Z_t \Vert_r + \Vert Z'_t \Vert_r))
            \Psi(\tilde X_t).
\label{pos}
\end{eqnarray}
The remainder of the uniqueness part is mainly devoted to proving this
result.

\subsection{An upper bound}
\label{parupperbound}
The first step towards the nonnegativity of the sum \eqref{pos}
is to give a bound on the term involving $f$. This is the purpose of
Lemma \ref{upbddpsi} below.

Let $(x,x')$ be a point in $\omb \times \omb$ and $b \in \R^d$. For
notational convenience, we keep the same notation $\omb \times \omb$ for
the image of this compact set in local coordinates considered below, and
we write $f$ for $f(b,x,z)$ and $f'$ for $f(b,x',z')$ (note that the $b$
is the same in $f$ and $f'$).
Take a local chart $(\phi,\phi)$ in which $(x,x')$ has coordinates
$(\hat x, \hat x')$; if $({\partial_1}, \ldots, {\partial_{2n}})$
denotes the natural dual basis of these coordinates, then
\begin{equation}
 \left(
   \begin{array}{c}
     f  \\
     f'
   \end{array}
 \right)
 = \sum_{i=1}^n \left( f^i \partial_i + f'^i \partial_{i+n} \right).
\label{coordf1}
\end{equation}
Now $v=(v_1, \ldots, v_{2n})=(\hat x - \hat x', \hat x')$ are new
coordinates in which the diagonal $\Delta$ is represented by the equation
$\{ v_1= \cdots = v_n =0 \}$; moreover, if
$(\varphi_1, \ldots, \varphi_{2n})$ is the natural dual basis of the
$v$-coordinates, then $\varphi_i=\partial_i$ and
$\varphi_{i+n} = \partial_i + \partial_{i+n}$ for $i=1, \ldots, n$;
thus \eqref{coordf1} becomes
$$
 \left(
   \begin{array}{c}
     f  \\
     f'
   \end{array}
 \right)
 = \sum_{i=1}^n \left( (f^i -f'^i) \varphi_i + f'^i \varphi_{i+n} \right)
$$
and
\begin{equation}
D \Psi \cdot
 \left(
   \begin{array}{c}
     f  \\
     f'
   \end{array}
 \right)
 = \sum_{i=1}^n \left( \frac{\partial \Psi}{\partial v_i} (f^i-f'^i)
                 + \frac{\partial \Psi}{\partial v_{i+n}} f'^i \right).
\label{majdpsi2}
\end{equation}
Finally, we call $p(V)$ the projection of a vector $V$ onto $\Delta$ :
if $V=(v_1, \ldots, v_{2n})$, then
$p(V)=(0, \ldots, 0, v_{n+1}, \ldots, v_{2n})$.

\begin{slemma}
\label{upbddpsi}
Suppose that $\Psi(x,x') \approx \delta(x,x')^\nu$ on $\omb \times \omb$
where $\nu$ is an even positive integer (since $\Psi$ is smooth).
Then there is $C>0$ such that,
for all $b$ in $\R^d$, $x,x'$ in $\omb$ and $z,z'$ in $\R^{n \times d_W}$

\begin{equation}
 \left\vert D\Psi \cdot
                                   \left(
                                 \begin{array}{c}
                                     f  \\
                                     f'
                                  \end{array}
                                   \right)
 \right\vert \le
 C  \delta(x,x')^{\nu-1} \left( \delta(x,x') \vert f' \vert
           + \vert f - f' \vert \right).
\label{majdpsi}
\end{equation}
\end{slemma}

\begin{demo}
First remark that on $\omb \times \omb$, we have
\begin{equation}
\vert V -p(V) \vert \approx \delta(x,x').
\label{equivalbis}
\end{equation}
Indeed, let $(x,x') \in M \times M$; in $v$-coordinates, it is
represented by $V=(\hat x - \hat x', \hat x')$. Then
$\vert V-p(V) \vert = \vert (\hat x - \hat x', 0) \vert$, the Euclidean
distance between $\hat x$ and $\hat x'$; as $x$ and $x'$
belong to the compact set $\omb$, the Euclidean and Riemannian distances
are equivalent and \eqref{equivalbis} follows.

Write $D_i$ for $\frac{\partial}{\partial v_i}$.
Since $\Psi(x,x') \approx \delta(x,x')^\nu$, a Taylor expansion gives
$$ \Psi(V) = \frac{1}{\nu !}
   \sum_{1 \le i_1, \dots, i_\nu \le n} D_{i_1 \dots i_\nu} \Psi (p(V))
   (V-p(V))_{i_1} \dots (V-p(V))_{i_\nu}
                     + {\it O}(\vert V-p(V) \vert^{\nu+1}).$$
Using another Taylor expansion with $D_i \Psi$ gives
\begin{eqnarray*}
\left\vert D_i \Psi (V) \right\vert
& \le &
    C \vert V -p(V) \vert^{\nu-1} \ \ \hbox{if} \ i \le n \\
\left\vert D_i \Psi (V) \right\vert
& \le &
    C \vert V -p(V) \vert^\nu \ \ \hbox{if} \ i > n \
         \hbox{(i.e. in the direction of the diagonal)}
\end{eqnarray*}
with a uniform $C$ on the compact set $\omb \times \omb$.
As a consequence of \eqref{majdpsi2}, we obtain
$$
 \left\vert D\Psi \cdot
                                   \left(
                                 \begin{array}{c}
                                     f  \\
                                     f'
                                  \end{array}
                                   \right)
 \right\vert \le
 C \vert V-p(V) \vert^{\nu-1} \left( \vert V-p(V) \vert
           \cdot \vert f' \vert + \vert f-f' \vert \right).
$$
Using once again \eqref{equivalbis} completes the proof.
\end{demo}

We now give two functions which verify the hypothesis of Lemma
\eqref{upbddpsi}. In these two examples, the Levi-Civita connection is
used.

\bigskip
\begin{expl}{\bf 1.}
Consider the mapping $\Psi (x,x') = \delta^2(x,x')$; it is smooth, not
convex in general, but this is true if the sectional curvatures are
bounded above by $0$.
\end{expl}

\medskip
\begin{expl}{\bf 2.}
Take for $\omb$ a regular geodesic
ball $\cB$ centered at $o \in M$, with the
sectional curvatures bounded from above by a constant $K>0$;
then we can find a nonnegative convex function $\Psi$ on $\cB \times \cB$
which vanishes only on the diagonal (see \cite{kend91})~:
$$ \Psi(x,x')= \left(
  \frac{1- \cos (\sqrt K \delta(x,x'))}
       {\cos (\sqrt K  \delta(x,o)) \cos (\sqrt K \delta(x',o)) - h^2 }
    \right)^p,$$
where $h>0$ is small and $p \ge 2$ is an integer large enough (so that
$\Psi$ is smooth).
\end{expl}

\subsection{The case $f$ independent of $z$}
In all this paragraph, $f$ does not depend on $z$, i.e.
$f(b,x,z)=f(b,x)$.

\begin{sprop}
\label{p1}
If the drift $f$ doesn't depend on $z$, then the process
$(e^{\lambda t} \Psi(\tilde X_t))_t$
is a submartingale for $\lambda$ positive large enough;
this implies that two $\omb$-valued
continuous  semimartingales $(X_t)$ and $(X'_t)$ verifying the same
equation $(M+D)$,  with the same terminal value, are indistinguishable.
\end{sprop}

\begin{demo}
For $x,x' \in \omb$, $\vert f(b,x) - f(b,x') \vert \le L' \delta(x,x')$;
thus using \eqref{majdpsi} and the boundedness of $f$

$$\left\vert D\Psi(x,x') \cdot
                                   \left(
                                 \begin{array}{c}
                                     f(b,x)  \\
                                     f(b,x')
                                  \end{array}
                                   \right)
 \right\vert \le C \delta(x,x')^\nu \le \tilde C \Psi(x,x').$$

Then for $\lambda \ge \tilde C$, $\mu=0$, the sum \eqref{pos} is
nonnegative (note that the convexity of $\Psi$ gives the nonnegativity
of the term involving $\Hess \Psi$). Moreover, the local martingale in
equation \eqref{ito3} is in fact a martingale since $D \Psi$ is bounded
on $\omb \times \omb$, so the process $(e^{\lambda t} \Psi(\tilde X_t))_t$
is indeed a submartingale. As it is nonnegative and has terminal value $0$,
it vanishes identically; so $\Psi(\tilde X_t)=0$ for all $t$. Finally, the
definition of $\Psi$ leads to $X_t=X'_t$ for all $t$ and the proof is
completed since we consider continuous processes.
\end{demo}

{\bf Remark. } Of course, $X=X'$ implies that for any $t$, $Z_t=Z'_t$
a.s.

As an immediate corollary, we are now able to give the uniqueness
property.
\begin{stheorem}
\label{unicit1}
Suppose that $\omb$ and $\Psi$ verify the properties of the
introduction (see paragraph \ref{par1.3}) and moreover that the drift $f$
depends only on $(b,x)$ and verifies \eqref{lip} and \eqref{upperboundf}.
Then, for a
given terminal value $U$ in $\omb$, there is at most one $\omb$-valued
solution to the equation $(M+D)$.
\end{stheorem}

{\bf Example. } Suppose that the Levi-Civita connection is used. Then
Theorem \ref{unicit1} implies uniqueness on any compact set of a
Cartan-Hadamard manifold, and on any regular geodesic ball (with $K>0$);
indeed it suffices to consider respectively the functions $\delta^2$ and
$\Psi$ defined after Lemma \ref{upbddpsi}.

\subsection{The general case in nonpositive curvatures}
In this subsection, the drift $f$ depends also on $z$,
$\Psi=\half \delta^2$ and the connection used is Levi-Civita's one;
moreover $M$ is supposed to be a Cartan-Hadamard manifold (i.e.
simply connected with nonpositive sectional curvatures); remark then
that any closed geodesic ball is regular. By achieving explicit
calculations we are going to derive the uniqueness property for any
compact set. \\
The problem is to show that the process
$(\exp (A_t) \Psi (\tilde X_t))_t$
is a submartingale. But to define such a process, we need to
consider solutions in some class which we now define.

\begin{sdefin}
\label{d1}
If  $\al$ is a positive constant, let $(\cal E_\al)$ be the set of
$\omb$-valued solutions of $(M+D)$ satisfying
\begin{equation}
  \esp \exp \left( \al \int_0^T \Vert Z_s \Vert_r^2 ds \right)
                 < \infty.
\label{definteg}
\end{equation}
\end{sdefin}

Actually, we now verify that for $\al$ small, $(\cal E_\al)$ contains
any solution of equation $(M+D)$.
The first step is the following lemma, which generalizes
Proposition 2.1.2 of \cite{mpl}.

\begin{slemma}
\label{Emu}
Suppose that we are given a positive
constant $\al$ and a $C^2$ function $\phi$ on $\omb$ satisfying
$C_{min} \le \phi(x) \le C_{max}$ for some positive $C_{min}$ and
$C_{max}$. Suppose moreover that
$\Hess \phi + 2 \al \phi \le 0$ on $\omb$; this means that
\begin{equation}
  \Hess \phi(x)<u,u>+2 \al \phi(x) \vert u \vert_r^2 \le 0.
\label{min1}
\end{equation}
Then, for every $\ep>0$, any $\omb$-valued solution of $(M+D)$ belongs
to $(\cal E_{\al - \ep})$.
\end{slemma}

\begin{demo}
Define
$$ S_t = \phi (X_t) \exp \left(
   \al \int_0^t \Vert Z_s \Vert_r^2 ds
 - \frac{1}{C_{min}} \int_0^t \vert D \phi(X_s) \cdot f(B_s^y,X_s,Z_s)
    \vert ds \right).$$
Denote by $e(t)$ the exponential term above. It follows from It\^o's
formula that
\begin{eqnarray*}
d S_t = e(t)
&\Bigg(&  D \phi(X_t) (Z_t d W_t)                                  \\
&      & + \left( \half \sum_{i=1}^{d_W} {}^t [{}^t Z_t]^i \Hess \phi (X_t)
    [{}^t Z_t]^i + \al \phi(X_t) \Vert Z_t \Vert_r^2 \right) dt       \\
&      & + \left( - \frac{\phi(X_t)}{C_{min}} \vert D \phi(X_t) \cdot
          f(B_t^y,X_t,Z_t) \vert
                + D \phi(X_t) \cdot f(B_t^y,X_t,Z_t) \right) dt \Bigg).
\end{eqnarray*}
Thus condition \eqref{min1} ensures that $S_t$ is a local
supermartingale; since it is nonnegative, we get
$\esp S_T \le \esp S_0$. By using the lower and upper bounds on
$\phi$, we deduce that
\begin{equation}
  \esp \exp \left(
       \al \int_0^T \Vert Z_s \Vert_r^2 ds - \frac{1}{C_{min}} \int_0^T
       \vert D \phi(X_s) \cdot f(B_s^y,X_s,Z_s) \vert ds
\right)
                 \le \frac{C_{max}}{C_{min}}.
\label{integ}
\end{equation}
By conditions \eqref{upperboundf} and \eqref{lip2}, there is a
$\tilde C>0$ such that
\begin{eqnarray*}
\frac{1}{C_{min}} \vert D \phi(X_s) \cdot f(B_s^y,X_s,Z_s) \vert
& \le & \tilde C (1+ \Vert Z_s \Vert_r)            \\
& \le & \ep \Vert Z_s \Vert_r^2 + C_\ep;
\end{eqnarray*}
in particular, we have \eqref{definteg} for $\al - \ep$.
\end{demo}

\begin{slemma}
\label{l2}
Let $\overline{B(o,\rho)}$ be a (regular) geodesic ball of center $o$ and
radius $\rho$. Then
for $\al>0$ small enough, every $\overline{B(o,\rho)}$-valued solution of
$(M+D)$ belongs to the set $(\cal E_\al)$ defined previously.
\end{slemma}

\begin{demo}
Let $\phi(x)=\cos \left( \frac{\pi}{3\rho} \delta(o,x) \right)$ for $x$
in the geodesic ball $\overline{B(o,\rho)}$.
We want to prove that this function $\phi$ verifies the hypothesis of
Lemma \eqref{Emu}.
Firstly, $\phi$ is obviously a smooth function and there is a $c>0$ such
that
$$\forall x \in \overline{B(o,\rho)}, \ c \le \phi(x) \le 1.$$
Moreover, we have
\begin{eqnarray*}
\Hess \phi (x) <u,u> =
& - & \cos \left( \frac{\pi}{3 \rho} \delta(o,x) \right)
    \left( \frac{\pi}{3 \rho} \right)^2 ( \delta'_2(o,x)<u>)^2       \\
& - & \sin \left( \frac{\pi}{3 \rho} \delta(o,x) \right)
     \left( \frac{\pi}{3 \rho} \right)
          {\rm Hess}_{22} \  \delta (o,x)<u,u>.
\end{eqnarray*}
If we get back to the notations of Subsection \ref{estimdist}, we
have, as a consequence of \eqref{der1} and \eqref{der2}, the estimates
(for $x \not= o$)~:
$$ \vert \delta'_2(o,x)<u> \vert = \vert v \vert_r$$
and
$${\rm Hess}_{22} \  \delta (o,x)<u,u> \ge
      \frac{\vert w \vert_r^2}{\delta(o,x)}.$$
Then
\begin{eqnarray*}
\Hess \phi (x) <u,u>
& \le & - \phi(x) \left( \frac{\pi}{3 \rho} \right)^2 \vert v \vert_r^2
        - \frac{ \sin \left( \frac{\pi}{3 \rho} \delta(o,x) \right)}
               { \frac{\pi}{3 \rho} \delta(o,x)}
         \left( \frac{\pi}{3 \rho} \right)^2 \vert w \vert_r^2          \\
& \le & - \left( \frac{\pi}{3 \rho} \right)^2 \phi(x) \vert u \vert_r^2.
\end{eqnarray*}
So there is a $\al>0$ (depending on the radius $\rho$)
such that \eqref{min1} holds. It suffices to apply Lemma \ref{Emu}
to conclude.
\end{demo}

Obviously, this is also true for any compact set $\omb$.
We can easily derive, from \eqref{integ} and the preceding lemma, the
following result, which will be useful in later calculations.
\begin{scorol}
\label{unifinteg}
There is an $\al>0$ and a finite positive constant $C_u$ (both depending
only on $\omb$ and the constants $L$ and $L_2$ in \eqref{lip} and
\eqref{upperboundf})
such that for any $\omb$-valued solution $(X,Z)$ of
equation $(M+D)$,
$$
  \esp \exp \left( \al \int_0^T \Vert Z_s \Vert_r^2 ds \right)
                 \le C_u.
$$
\end{scorol}

We can now state the uniqueness result (Theorem \ref{unicit2}).

\begin{slemma}
\label{l3}
For two solutions $(X,Z)$, $(X',Z')$ of $(M+D)$, verifying
$\esp (A_T) < \infty$
for every $\mu>0$, the expression \eqref{pos} is nonnegative for
$\lambda$ and $\mu$ large enough.
\end{slemma}

\begin{demo}
Using
\eqref{min} we have for
$\tilde z= \left(
          \begin{array}{c}
             z  \\
             z'
          \end{array} \right)$
\begin{equation}
  \sum_{i=1}^{d_W} {}^t [{}^t \tilde z]^i
        \Hess \Psi (\tilde x) [{}^t \tilde z]^i   \ge
  \sum_{i=1}^{d_W} \Big\vert \overset{x'}{\underset{x}{\Vert}} [{}^t z]^i -
       [{}^t z']^i \Big\vert_r ^2  =
 \Big\Vert \overset{x'}{\underset{x}{\Vert}} z -z' \Big\Vert_r ^2.
\label{minhess}
\end{equation}

Moreover, we also have using \eqref{lip} and \eqref{upperboundf},
together with \eqref{majdpsi} (or with \eqref{tp1} and \eqref{der1})
\begin{eqnarray}
 \left\vert D\Psi (\tix)
                                   \left(
                                 \begin{array}{c}
                                     f(b,x,z)    \\
                                     f(b,x',z')
                                  \end{array}
                                   \right)
 \right\vert
& \le &
 C_1 \delta(x,x') \left(
 \delta(x,x')(1+ \Vert z \Vert_r + \Vert z' \Vert_r)
         +\left\Vert \overset{x'}{\underset{x}{\Vert}} z -z'
         \right\Vert_r \right)                           \nonumber   \\
& \le & C \delta^2(x,x')(1+\Vert z \Vert_r +
                                   \Vert z' \Vert_r)
       + \frac{1}{4} \left\Vert
          \overset{x'}{\underset{x}{\Vert}} z -z' \right\Vert_r^2.
\label{majdpsi3}
\end{eqnarray}
Then \eqref{pos} is greater than the following sum
$$-C \delta^2(X_t,X_t')(1+\Vert Z_t \Vert_r
                       + \Vert Z'_t \Vert_r)
   + \half \delta^2(X_t,X_t')(\lambda +
              \mu(\Vert Z_t \Vert_r + \Vert Z'_t \Vert_r)).$$

Taking $\lambda$ and $\mu$ greater than $2C$ makes obviously this sum
(and the expression \eqref{pos}) nonnegative.
\end{demo}

\begin{stheorem}
\label{unicit2}
Suppose that $M$ is a Cartan-Hadamard manifold and that the drift $f$
verifies condition \eqref{lip} and \eqref{upperboundf}. Then, for a
given terminal value $U$ in the compact $\omb$, there is at most one
$\omb$-valued solution to the equation $(M+D)$ (i.e. for two solutions
$(X,Z)$ and $(X',Z')$, the (continuous) processes $X$ and $X'$ are
indistinguishable).
\end{stheorem}

\begin{demo}
 Every compact $\omb$ is included in a closed geodesic ball. Hence any
$\omb$-valued solution of $(M+D)$ is in $({\cal E}_\al)$ for a $\al>0$
from Lemma \ref{l2}; this easily gives
the integrability of $\exp(A_T)$ for every $\mu>0$. Consequently, it 
suffices to apply Lemma \ref{l3} and conclude as in the proof of 
Proposition \ref{p1}.
\end{demo}

{\bf Remark.} Note that the only hypothesis required for uniqueness
in this case is the compactness of $\omb$.

\section{Existence results}
\label{par4}

In this section we are given an $\omb$-valued random variable $U$
and we want to construct a pair of processes $(X,Z)$, satisfying
equation $(M+D)$, with $X$ in $\omb$ and terminal value $U$.
We limit ourselves to the case of a
Wiener probability space and we recall that on $\omb$, if $(X,Z)$
and $(X',Z')$ are two solutions of the equation $(M+D)$ and
$\tilde X =(X,X')$, then for $\lambda>0$ and $\mu>0$ large enough
the following processes
$$
\left( e^{\lambda t} \Psi(\tilde X_t) \right)_{t \in [0;T]}
\ \ \hbox{ if $f=f(B^y,X)$ }
$$
and, when $M$ is Cartan-Hadamard,
$$
\left( e^{\lambda t +
          \mu \int_0^t (\Vert Z_s \Vert_r + \Vert Z'_s \Vert_r) ds}
           \delta^2(\tilde X_t) \right)_{t \in [0;T]}
\ \ \hbox{ if $f=f(B^y,X,Z)$ }
$$
are nonnegative submartingales.

The strategy of the proof can be described as follows~:

\noindent
1. Simplify the problem by considering only terminal values which
can be expressed as functions of the diffusion $B^y$ at time $T$,
i.e. $U=F(B_T^y)$ (Subsection \ref{par4.1}). This step needs to pass
through the limit in equation $(M+D)$; when $f$ is independent of $z$,
it is a corollary of a well-known result, but in the other case, more
technical calculations, using the uniqueness part, are involved.

\noindent
2. Solve a Pardoux-Peng BSDE with parameter to construct a pair of
processes in $\R^n \times \R^{n d_w}$ which is close to being a solution
of $(M+D)$ with $X_T=U$ (Subsection \ref{par4.2}).

\noindent
3. Show that under an additional condition on $f$ the solution of
the preceding BSDE is a solution of the BSDE $(M+D)$ on a small time
interval (Subsections \ref{par4.3} and \ref{par4.4}). Note that the main
argument in the existence proof is certainly Proposition \ref{zbounded},
where we give an a.s. upper bound on the process $(Z_t)$.

\noindent
4. Use the convex function $\Psi$ to show that we have a solution of
$(M+D)$ on the whole time interval $[0;T]$  (Subsection \ref{par4.5}).

In fact, for technical reasons we suppose in the two last steps that
$f$ is sufficiently regular; then the proof of the existence is
completed with the last subsection~:

\noindent
5. Solve BSDE $(M+D)$ for general $f$ using classical approximation
methods (Subsection \ref{par4.6}).

Note that we usually work within local coordinates in $\R^n$, i.e. we
consider that $\omb \subset O \subset \R^n$.

\subsection{Reduction of the problem}
\label{par4.1}
Let $C^\infty_c (\R^d, \omb)$ denote $\omb$-valued functions on $\R^d$
which are constant off a compact set. In this paragraph, it is shown
that it suffices to check the existence result for $U=F(B_T^y)$ with
$F \in C^\infty_c (\R^d, \omb)$.

\bigskip

The space of all functionals
$\{ G(W_{t_1}, W_{t_2}, \ldots , W_{t_q}),
      0<t_1 < \ldots < t_q \le T ; \ G \in
                      C^\infty_c (\R^{q d_W}, \omb) \}$ is dense in
$L^2(\cF_T ; \omb)$ (remark that
$L^2(\cF_T ; \omb) = L^\infty(\cF_T ; \omb)$ since $\omb$ is compact)
endowed with the distance

$$ \delta^{(1)} (U_1, U_2) =
          \sqrt {\esp \left( \delta^2 (U_1, U_2) \right)}.$$
The same is true by replacing functions of the Brownian Motion at
discrete times by functions of the diffusion $B^y$ at discrete times.
Indeed, if we add the $d_W$ components $W^1_t, \ldots, W^{d_W}_t$
to the $d$ components of the diffusion
$B_t^y=(B_t^{1,y}, \ldots, B_t^{d,y})$, it is easy to
conclude that the diffusion obtained in this way generate the same
filtration as $(W_t)_{0 \le t \le T}$.

Let $U^l \in L^2({\cal F}_T; \omb)$ and $(X^l,Z^l)$ the $\omb$-valued
solution of $(M+D)$ with $X_T^l=U^l$. We are going
to show that if $U^l \rightarrow U$ in $L^2({\cal F}_T; \omb)$, then
$(X^l)$ tends, for the distance
$$ \delta^{(2)} \left( (X^l_t), (X_t) \right) =
   \sqrt {\esp \left( \sup_{t \in [0,T]} \delta^2 (X^l_t, X_t)
                                                 \right)},$$
to a process $X$ ending at $U$ and that there is a process $Z$ such that
$(X,Z)$ solves the BSDE $(M+D)$.

\begin{sdefin}
Let ${\cal T}$ be the set of all terminal values $U \in \omb$ of
processes $X$ solutions of $(M+D)$ (i.e. such that there is a process $Z$
with $(X,Z)$ solution of $(M+D)$) and ${\cal S}$ be the set of all
these $\omb$-valued processes.
According to the uniqueness part, to every $U$ in ${\cal T}$ corresponds
a unique process $(X_t)_t$ in ${\cal S}$ such that $X_T = U$. Hence we
can define a mapping $c$ with $c(U)$ being this process~:

$$\begin{array}{cccc}
c : & {\cal T} & \rightarrow & {\cal S} \\
    & U        & \mapsto     & (X_t)_t.
\end{array}$$

We endow the space ${\cal T}$ with $\delta^{(1)}$ and ${\cal S}$ with
$\delta^{(2)}$, the distances just defined.
\end{sdefin}

It is obvious that $c$ is one-to-one and onto, and
that $c^{-1}$ is uniformly continuous for the above distances.
Next we want to prove the uniform continuity of $c$.
This is the aim of the following lemma, whose proof is a (slightly)
modified version of the one of theorem (5.5) from \cite{kend90}.

\begin{slemma}
$c$ is uniformly continuous for the distances $\delta^{(1)}$
and $\delta^{(2)}$.
\label{unifcontin}
\end{slemma}

\begin{demo}
Suppose that we are given two terminal values $U_1$ and $U_2$
corresponding to two solutions $(X_t,Z_t)$ and $(X'_t,Z'_t)$.
Then by H\"older's inequality,
\begin{equation}
\delta^{(2)}((X_t),(X'_t)) \le \esp \left(
                     \sup_t \delta^p(X_t,X'_t) \right) ^\frac{1}{p},
\label{major1}
\end{equation}
and since $\Psi \approx \delta^p$
\begin{eqnarray}
\esp \left( \sup_t \delta^p(X_t,X'_t) \right)
& \le & C \esp \left( \sup_t
        \left( e^{\lambda t +
          \mu \int_0^t (\Vert Z_s \Vert_r + \Vert Z'_s \Vert_r) ds}
           \Psi(X_t,X'_t) \right)^2 \right)^\half       \nonumber  \\
& \le & C \esp \left( e^{2\lambda T +
          2\mu \int_0^T (\Vert Z_s \Vert_r + \Vert Z'_s \Vert_r) ds}
           \Psi^2(X_T,X'_T) \right)^\half               \nonumber  \\
& \le & C \esp \left( e^{ 4 \mu \int_0^T (\Vert Z_s \Vert_r
         + \Vert Z'_s \Vert_r) ds} \right) ^\frac{1}{4}
         \esp \left( \delta^{4p}(U_1,U_2) \right) ^\frac{1}{4}
                                                \nonumber       \\
& \le & C \esp \left( e^{ C_\eta T + \eta \int_0^T (\Vert Z_s \Vert_r ^2
         + \Vert Z'_s \Vert_r^2) ds} \right) ^\frac{1}{4}
         \esp \left( \delta^2(U_1,U_2) \right) ^\frac{1}{4}.
\label{major2}
\end{eqnarray}
The constant $C$ above is allowed to vary from one inequality to
another, but it depends only on $T$, $\omb$ and $\Psi$ (but not on
the processes $X$ and $X'$).
The second inequality is Doob's $L^2$ one applied to the submartingale
$(\exp(A_t) \Psi(X_t,X'_t))_t;$
the third one is Cauchy-Schwarz's one and the last one uses the
classical inequality $x \le \eta x^2 + C_\eta$
and the boundedness of $\delta$ which implies
$\delta^{4p} \le \tilde C \delta^2$.

Inequalities \eqref{major1} and \eqref{major2} together give~:
\begin{equation}
\delta^{(2)}((X_t),(X'_t)) \le C \esp \left( e^{ C_\eta T
         + \eta \int_0^T (\Vert Z_s \Vert_r^2
         + \Vert Z'_s \Vert_r^2) ds} \right) ^\frac{1}{4p}
      \delta^{(1)}(U_1,U_2)^\frac{1}{2p}.
\label{major3}
\end{equation}
Now we specialize to each case~:
if $f$ does not depend on $z$, take $\mu=0$; then
$\left( e^{\lambda t} \Psi(X_t,X'_t) \right)_t$ is indeed a
submartingale, and with $\eta=C_\eta=0$, the inequality
\eqref{major3} gives the uniform continuity required.

In the other case,
$(\exp(A_t) \Psi(X_t,X'_t))_t$
is yet a submartingale with $\Psi=\half \delta^2$; moreover,
taking in \eqref{major3} $\eta=\al$ (with $\al$ as in Corollary
\ref{unifinteg}) leads to the conclusion again.
\end{demo}

Using the completeness of the space of all $\omb$-valued processes
endowed
with the  distance $\delta^{(2)}$, we get the following result as
an easy consequence of the preceding lemma.

\begin{sprop}
Let $(U^l)_l$ be a sequence in ${\cal T}$ converging to
$U \in L^2(\cF_T ; \omb)$ for the distance $\delta^{(1)}$ and
$X^l=c(U^l)$.
Then there is a (continuous) process $(X_t)_{t \in [0;T]}$ such that
$\delta^{(2)}(X^l,X)\overset{l}{\rightarrow} 0$.
In particular, we have $X_T = U$ a.s.
\label{XlcvX}
\end{sprop}

The result of this part is completed by proving the next proposition.

\begin{sprop}
There is a process $Z$ such that the pair $(X,Z)$, with $X$ defined in
Proposition \ref{XlcvX}, is a solution of $(M+D)$ with $X_T=U$, limit
of the sequence $(U^l)_l$.
\label{passlim}
\end{sprop}

\begin{demo}
First, it is clear from the definition of $\delta^{(2)}$ that $X$
is $\omb$-valued.
We give then two proofs, one for the $z$-independent case and one for
the general case in nonpositive curvatures.

In the $z$-independent case, we use a proof very similar to the one of
Theorem (4.43) of \cite{scm} based on the characterization of
Proposition \ref{charact}~: first we
localize the processes (in order to work in the $O_p$ defined in
Proposition \ref{charact}) and then it suffices to notice that the
submartingale property passes through the limit in $l$.

\bigskip
In the other case, it is not possible to apply the preceding proof, as
the drift depends also on $z$. In fact, we are going to find explicitly
a process $Z$ such that $(X,Z)$ is a solution of $(M+D)$.
We recall that for each $l$, we have in $\R^n$ (remember that
$[Z^l_t]^j$ is the $j^{th}$ row of the matrix $Z^l_t$)
$$\left\{       \begin{array}{l}
  d X_t^l = Z^l_t d W_t                    + \left(- \half
\Gamma_{jk}(X^l_t)             (\left[ Z^l_t \right]^k \vert
         \left[ Z^l_t \right]^j) + f(B_t^y,X_t^l, Z_t^l)  \right) d t \\
	  X_T^l=U^l.                                         \\
      \end{array}
    \right.   $$

{ \it First step : } Find a process $Z$, limit in the $L^2$ sense of the
processes $Z^l$.

Let
$$\tilde X^{l,m}=(X^l,X^m) \ \ \hbox{ and } \ \
  \tilde Z^{l,m} = \left(
                 \begin{array}{c}
                    Z^l   \\
	            Z^m
                  \end{array}
                \right);$$
apply It\^o's formula \eqref{ito1} to $\Psi(\tilde X^{l,m})$; then use
\eqref{minhess} to bound below the Hessian term and \eqref{majdpsi3} to
bound above the term involving $f$. Taking the expectation gives

\begin{eqnarray*}
\frac{1}{4} \esp \int_0^T \left\Vert
  \overset{X_s^l}{\underset{X_s^m}{\Vert}} Z_s^m - Z_s^l
                                         \right\Vert_r^2 ds
&\le& C \esp \left( \sup_s \delta^2(X_s^l,X_s^m) \right)              \\
&   & + C \esp \int_0^T \delta^2(X_s^l,X_s^m)
             (\Vert Z_s^l \Vert_r + \Vert Z_s^m \Vert_r) ds.
\end{eqnarray*}
We know that the first expectation on the right tends to zero
as $l$ and $m$ tend to $+ \infty$. Moreover, if  $I_1$ denotes the last
integral, then using a $\al$ as in Corollary \ref{unifinteg} and the
fact that $\delta^4$ is bounded above by $C \delta^2$ on the compact set
$\omb \times \omb$, we get
\begin{eqnarray}
I_1
& \le & \sqrt 2 \esp \left( \sup_s \delta^4(X_s^l,X_s^m) \right)^\half
        \esp \left( \int_0^T
            (\Vert Z_s^l \Vert_r^2 +
               \Vert Z_s^m \Vert_r^2) ds \right)^\half    \nonumber   \\
& \le & \sqrt{\frac{2}{\al}} \esp \left(
      \sup_s \delta^4(X_s^l,X_s^m) \right)^\half
   \esp \left( e^{\al \int_0^T (\Vert Z_s^l \Vert_r^2
        + \Vert Z_s^m \Vert_r^2) ds} \right)^\half        \nonumber   \\
& \le & \frac{C}{\sqrt \al} \esp \left(
      \sup_s \delta^2(X_s^l,X_s^m) \right)^\half.
\label{majI1}
\end{eqnarray}

As $C$ does not depend on $l,m$,
\begin{equation}
\esp \int_0^T \left\Vert
  \overset{X_s^l}{\underset{X_s^m}{\Vert}} Z_s^m - Z_s^l \right\Vert_r^2 ds
       \overset{l,m \rightarrow \infty}{\longrightarrow} 0.
\label{cvZnnormr}
\end{equation}

Now we use inequality \eqref{tp2} to have a bound on the Euclidean
norm~:
\begin{eqnarray*}
 \esp \int_0^T \Vert Z_s^m - Z_s^l \Vert ^2 ds
 \le  C \bigg(
&   & \esp \int_0^T \left\Vert \overset{X_s^l}{\underset{X_s^m}{\Vert}}
         Z_s^m - Z_s^l \right\Vert_r^2 ds            \\
& + & \esp \int_0^T \delta^2(X_s^l,X_s^m) (\Vert Z_s^l \Vert_r^2
                   + \Vert Z_s^m \Vert_r^2 ) ds \bigg).
\end{eqnarray*}

The first term on the right tends to $0$ by \eqref{cvZnnormr}, and
an argument similar to \eqref{majI1} would show that it also holds
for the second term. Hence

$$
 \esp \int_0^T \Vert Z_s^m - Z_s^l \Vert ^2 ds
    \overset{l,m \rightarrow \infty}{\longrightarrow} 0.
$$
Now by completeness of the space $L^2([0;T] \times \Omega)$, we have
the required result :
$$
\exists (Z_t)_t \in L^2([0;T] \times \Omega) : (Z_t^l)
     \overset{L^2}{\longrightarrow} (Z_t).
$$

\bigskip
{ \it Second step : } $(X,Z)$ is indeed a solution of
equation $(M+D)$ with terminal value $U$.

In view of this, let us show that the following expectation tends
to zero as $l$ tends to $\infty$ :
\begin{eqnarray*}
\esp \bigg\vert
U & - & \int_t^T Z_s d W_s - \int_t^T \left( - \half
   \Gamma_{jk}(X_s) (\left[ Z_s \right]^k \vert
         \left[ Z_s \right]^j) + f(B_s^y,X_s, Z_s)  \right) d s        \\
  & - & U^l + \int_t^T Z^l_s d W_s + \int_t^T \left( - \half
   \Gamma_{jk}(X^l_s) (\left[ Z^l_s \right]^k \vert
         \left[ Z^l_s \right]^j) + f(B_s^y,X_s^l, Z_s^l)  \right) d s
        \bigg\vert.
\end{eqnarray*}
Obviously, this expectation is bounded above by
\begin{eqnarray*}
 \esp \left( \vert U-U^l \vert^2 \right)^\half
& + &\esp\left( \int_0^T \Vert Z_s^l-Z_s \Vert^2 ds \right)^\half\\
& + & \esp \left( \int_0^T \vert \Gamma_{jk}(X_s) - \Gamma_{jk}(X_s^l)
     \vert \cdot \vert (\left[ Z_s \right]^k \vert
     \left[ Z_s \right]^j) \vert ds \right)                    \\
& + & \esp \left( \int_0^T \vert \Gamma_{jk}(X_s^l) \vert
     \left\vert  (\left[ Z_s \right]^k \vert \left[ Z_s \right]^j)
        -(\left[ Z_s^l \right]^k \vert \left[ Z_s^l \right]^j)
       \right\vert ds \right)                                  \\
& + & \esp \left( \int_0^T \vert f(B_s^y,X_s, Z_s)
                          - f(B_s^y,X_s^l, Z_s^l)  \vert ds \right).
\end{eqnarray*}
We know that the first two expectations tend to zero;
the third
term tends to zero by dominated convergence (at least for a subsequence
of $(X^l)$, but it doesn't matter since $(X^l)$ is Cauchy).
Let $E_1$ denote the next expectation; then we can write
\begin{eqnarray*}
E_1
& \le & C \esp \left( \int_0^T \Vert Z_s^l - Z_s \Vert
    (\Vert Z_s^l \Vert + \Vert Z_s \Vert) ds \right)           \\
& \le & \sqrt 2 C \esp \left( \int_0^T \Vert Z_s^l - Z_s \Vert ^2 ds
        \right) ^\half
        \esp \left( \int_0^T (\Vert Z_s^l \Vert^2 + \Vert Z_s \Vert^2)
             ds \right) ^\half.
\end{eqnarray*}
The first integral tends to zero and the second is bounded because $Z^l$
converges in $L^2$; hence $E_1$ tends to zero.

Finally, let $E_2$ denote the last integral and use inequality
\eqref{lip2} to obtain
$$
E_2
\le \esp \left( \int_0^T \left( \delta(X_s^l,X_s)(1+ \Vert Z_s^l \Vert
      + \Vert Z_s \Vert) + \Vert Z_s^l - Z_s \Vert \right) ds \right);
$$
using an argument similar to \eqref{majI1}, this quantity tends also to
zero as $l$ tends to $\infty$.

Hence the limit in $L^1$ of $X_t^l$ is
$$U - \int_t^T Z_s d W_s - \int_t^T \left( - \half
   \Gamma_{jk}(X_s) (\left[ Z_s \right]^k \vert
         \left[ Z_s \right]^j) + f(B_s^y,X_s, Z_s)  \right) d s;$$
We know that it is also $X_t$, so by continuity~:
$$ a.s., \ \forall t, \
X_t=U - \int_t^T Z_s d W_s - \int_t^T \left( - \half
   \Gamma_{jk}(X_s) (\left[ Z_s \right]^k \vert
         \left[ Z_s \right]^j) + f(B_s^y,X_s, Z_s)  \right) d s.$$
That finishes the proof of the proposition for a general $f$ in
nonpositive curvatures.
\end{demo}

As a consequence of Propositions \ref{XlcvX} and \ref{passlim}, it
suffices to work with a dense subset of $L^2(\cF_T ; \omb)$, i.e. to
consider terminal values $U$ that are written
$ U = G(B_{t_1}^y, B_{t_2}^y, \ldots , B_{t_q}^y)$ where
$G \in C^\infty_c (\R^{qd}, \omb)$ (in fact, we can obviously take
$t_q=T$).

A further step of simplification is possible (see \cite{kend90}):
conditioning by $\cF_{t_{q-1}}$ and working over the
time interval $[t_{q-1} ; T]$,
$B_{t_1}^y, \ldots , B_{t_{q-1}}^y$ can be treated as
constant. Then, if we know a solution $(X_t,Z_t)$ to $(M+D)$ on
$[t_{q-1} ; T]$ with $X_T=G(B_T^y)$, the problem is
to reach $X_{t_{q-1}}$. But this variable is in fact a measurable
function of $B_{t_1}^y, B_{t_2}^y, \ldots , B_{t_{q-1}}^y$.
A density argument enables us to suppose this function smooth and
constant off a compact set. Thus by induction the problem is solved
if we can find a solution to $(M+D)$ with terminal value $F(B_T^y)$
with $F \in C^\infty_c (\R^d, \omb)$.

So to prove the existence result for any $\omb$-valued terminal
variable $U$, it suffices to solve equation $(M+D)$ with a terminal
value $U$ that can be written $F(B_T^y)$ where $F \in C^\infty_c
(\R^d, \omb)$. This is the aim of the next paragraphs.

\subsection{Approximation by BSDEs with Lipschitz coefficients}
\label{par4.2}

We deal with a BSDE whose coefficients $\Gamma_{jk}(x)$ (we recall
that it is a vector in $\R^n$; see the introduction) and $f(b,x)$
(or $f(b,x,z)$) are defined only for $x$ in the open subset $O$ of
$\R^n$, and with a quadratic term in $Z_t$; but we would like to
apply the existence result of \cite{pp90} to BSDEs with Lipschitz
coefficients. Our purpose in this part is to define such a BSDE,
derived from the initial one and defined on all $\R^n$.

Firstly, we extend the definition of the BSDE $(M+D)$ to the whole
space $\R^n$~: let $\phi$  be a smooth function on $\R^n$ with
compact support in $O$ and such that $\phi=1$ on $\omb$.
We will explicit $\phi$ a little more later in Subsection \ref{par4.4}.
Then we extend the drift $f$ and Christoffel symbols to the whole
space $\R^n$ by letting , for $b \in \R^d$,
$\tilde f(b,x)= \phi(x) f(b,x)$ (or $\tilde f(b,x,z)= \phi(x)
f(b,x,z)$) and $\tilde \Gamma (x) = \phi(x) \Gamma(x)$.

\noi
The new BSDE defined on all $\R^n$ (we keep $X$ and $Z$ for the
notations) is :
$$\widetilde {(M+D)}
\left\{
      \begin{array}{l}
        d X_t = Z_t d W_t
                   + \left(- \half \tilde \Gamma_{jk}(X_t)
            (\left[ Z_t \right]^k \vert
     \left[ Z_t \right]^j) + \tilde f(B_t^y,X_t, Z_t)  \right) dt \\
	  X_T=U.                                                  \\
      \end{array}
    \right.   $$
{\bf Remark : } We write $\tilde f (B_t^y,X_t, Z_t)$ but when it is not
pointed out, it should also be interpreted as $\tilde f (B_t^y,X_t)$ as
well.

\bigskip
Now let $\ep \in ]0;1[$ and $s_\ep : \R_+ \rightarrow \R_+$ be a smooth
nondecreasing
function such that $s_\ep(t)=0$ iff $t \in [0; \frac{1}{\ep}]$, and
$s_\ep$ is linear for $t$ large enough. Then, for $z \in \R^{n d_W}$, we
define $h_\ep(z)= s_\ep(\Vert z \Vert)$ and
$\overline z = \frac{z}{1+ h_\ep(z)}$.
We have the following

\begin{slemma}
\label{coeflip}
The functions $h_\ep$ and $\overline z$ defined above satisfy the next
assertions.

\begin{item}{(i)}
If $\Vert z \Vert \le \frac{1}{\ep}$, then $h_\ep(z)=0 $ and
$\overline z=z$;
\end{item}

\begin{item}{(ii)}
$z \mapsto \overline z$ is smooth, bounded and Lipschitz on $\R^{n d_W}$;
\end{item}

\begin{item}{(iii)}
$(x,z) \mapsto g(x,z):=\tilde \Gamma_{jk} (x)
        ([\overline z]^k \vert [\overline z]^j)$ is smooth and Lipschitz
on $\R^n \times \R^{n d_W}$.
\end{item}

\begin{item}{(iv)}
$(b,x,z) \mapsto \tilde f (b,x,\overline z)$ is
bounded and Lipschitz on $\R^d \times \R^n \times \R^{n d_W}$.
\end{item}

As a consequence, the function
\begin{eqnarray*}
     \gamma : \R^d \times \R^n \times \R^{n d_W}
&\rightarrow& \R^n                               \\
(b,x,z)
& \mapsto &
\half \tilde \Gamma_{jk}(x)
  ([\overline z]^k \vert [\overline z]^j) - \tilde f (b,x,\overline z)
\end{eqnarray*}
is bounded and Lipschitz.

\end{slemma}

\begin{demo}

$(i)$ It is obvious on the definition.

$(ii)$ The smoothness and boundedness of $\overline z$ are clear.
Now we just prove that
$\Vert \overline z - \overline{z'} \Vert \le C \Vert z-z' \Vert$,
considering the 3 cases below.

$\bullet$ {\it $\Vert z \Vert \le \frac{1}{\ep}$ and
     $\Vert z' \Vert \le \frac{1}{\ep}$}

\noi Then $\Vert \overline z - \overline{z'} \Vert = \Vert z-z' \Vert$.

$\bullet$ {\it $\Vert z \Vert \ge \frac{1}{\ep}$ and
     $\Vert z' \Vert \ge \frac{1}{\ep}$}

\noi Then, as $\overline{z}$ is bounded and $h_\ep$ Lipschitz, we can
write
\begin{eqnarray*}
\left\Vert \overline z - \overline {z'} \right\Vert & \le &
    \left\Vert \frac{z-z'}{1+h_\ep(z)} \right\Vert
  + \Vert z' \Vert \  \left\vert \frac{1}{1+h_\ep(z)}
     - \frac{1}{1+h_\ep(z')} \right\vert                              \\
& \le &
   \Vert z-z' \Vert + \Vert \overline{z'} \Vert \ \left\vert
        \frac{h_\ep(z')-h_\ep(z)}{1+h_\ep(z)} \right\vert             \\
& \le & C \Vert z-z' \Vert.
\end{eqnarray*}

$\bullet$ {\it $\Vert z \Vert > \frac{1}{\ep}$ and
     $\Vert z' \Vert \le \frac{1}{\ep}$}

\noi Let $z'' \in [z;z']$ with $\Vert z'' \Vert = \frac{1}{\ep}$;
then, using the first two cases, we have
\begin{eqnarray*}
\left\Vert \overline z - \overline {z'} \right\Vert
& \le &
\left\Vert \overline z - \overline {z''} \right\Vert +
\left\Vert \overline {z'} - \overline {z''} \right\Vert
    \le C (\Vert z-z'' \Vert + \Vert z'-z'' \Vert) = C \Vert z-z' \Vert.
\end{eqnarray*}

$(iii)$ We have
\begin{eqnarray*}
\vert g(x,z)-g(x',z') \vert
& \le &
  \vert \tilde \Gamma_{jk} (x) - \tilde \Gamma_{jk} (x') \vert
  \cdot \vert ([\overline z]^k \vert [\overline z]^j) \vert    \\
&     & + \vert \tilde \Gamma_{jk} (x') \vert \cdot
  \vert ([\overline z]^k \vert [\overline z]^j)
     - ([\overline {z'}]^k \vert [\overline {z'}]^j) \vert.
\end{eqnarray*}
There is no problem with the first term on the right because
$\tilde \Gamma_{jk}$ is smooth on $\R^n$ with compact support
and $\overline z$ is bounded. Using (ii), it is easy
to deduce that the second term is bounded above by
$C \Vert z-z' \Vert$ and that finishes the proof.

$(iv)$ Recall that for $x, x' \in O$,
$$\vert f(b,x,\overline z) - f(b',x',\overline{z'}) \vert \le L'
    \left( (\vert b - b' \vert + \vert x - x' \vert) (1 + \Vert
                   \overline z \Vert + \Vert \overline{z'} \Vert)
           + \Vert \overline z - \overline{z'} \Vert \right).$$
If $\overline{O_s}$ denotes the compact support of $\phi$, this
inequality and \eqref{upperboundf} imply that
$f(b,x, \overline z)$ is bounded, uniformly in $x \in \overline{O_s}$,
$b \in \R^{d}$ and $z \in \R^{n d_W}$. Thus $\tilde f$ is bounded on
$\R^d \times \R^n \times \R^{n d_W}$ because for $x$ outside
$\overline{O_s}$, $\tilde f=0$.
Besides, we consider the 3 following cases to study the Lipschitz
property~:

$\bullet$ {\it $x, x' \notin \overline{O_s}$}

\noi Then $\tilde f(b,x,\overline z) - \tilde f(b',x',\overline{z'})=0$
and it is obvious.

$\bullet$ {\it $x \in \overline{O_s}, x' \notin \overline{O_s}$}

\noi Then, as $\phi$ is Lipschitz,
\begin{eqnarray*}
\vert \tilde f(b,x,\overline z) - \tilde f(b',x',\overline{z'}) \vert
& = & \vert \phi(x) f(b,x,\overline z) \vert                        \\
& = & \vert \phi(x)
              - \phi(x') \vert \cdot \vert f(b,x,\overline z) \vert \\
&\le& \alpha \vert x-x' \vert.
\end{eqnarray*}

$\bullet$ {\it $x, x' \in \overline{O_s}$}

\noi Then
\begin{eqnarray*}
\vert \tilde f(b,x,\overline z) - \tilde f(b',x',\overline{z'}) \vert
& = & \vert \phi(x) f(b,x,\overline z)
        - \phi(x') f(b',x',\overline {z'})\vert                      \\
&\le& \vert \phi(x)
          - \phi(x') \vert \cdot \vert f(b,x,\overline z) \vert      \\
&   & + \vert \phi(x') \vert \cdot \vert f(b,x,\overline z)
      -  f(b',x',\overline{z'}) \vert                                \\
&\le& \alpha \vert x-x' \vert
      + \beta_1 \Vert \overline z - \overline{z'} \Vert
      + \gamma \vert b-b' \vert           \\
&\le& \alpha \vert x-x' \vert + \beta_2 \Vert z - z' \Vert
      + \gamma \vert b-b' \vert
\end{eqnarray*}
since $z \mapsto \overline z$ is Lipschitz by $(ii)$.
The proof of $(iv)$ is completed.
\end{demo}

Now we can introduce a new BSDE
$$\widetilde {(M+D)_\ep}
\left\{
      \begin{array}{l}
        d X_t^\ep = Z^\ep_t d W_t
                - \gamma (B_t^y, X_t^\ep, Z_t^\ep) dt        \\
	  X_T^\ep=U.
      \end{array}
    \right.   $$
The interest of this new equation lies in the following result.

\begin{sprop}
\label{pardouxpeng}
The equation $\widetilde {(M+D)_\ep}$ has Lipschitz coefficients; it has
a unique
solution $(X^\ep_t,Z^\ep_t) \in \R^n \times \R^{n d_W}$ such that
$$\esp \left( \sup_{t \in [0;T]} \vert X_t^\ep \vert ^2 \right)<\infty
\ \hbox{ and } \
\esp \left( \int_0^T \Vert Z_t^\ep \Vert ^2 dt \right) <\infty.$$
\end{sprop}

\begin{demo}
The first part is a consequence of Lemma \ref{coeflip}; then existence
and uniqueness are classical results of \cite{pp90}.
\end{demo}

In Subsections \ref{par4.3}, \ref{par4.4} and \ref{par4.5} below, $f$ is
also
supposed to be a $C^3$ function such that
$(b,x,z) \mapsto \tilde f(b, x, \overline z)$ and $\gamma$ are $C^3$
functions
with all their partial derivatives of order 1, 2 and 3 bounded.

\subsection{Existence of a solution of $\widetilde {(M+D)}$ on a small
time interval}
\label{par4.3}

In this paragraph we consider a terminal value $U=F(B_T^y)$
($F \in C^\infty_c (\R^d, \omb)$) which lies in
$\omb$ and we show that for an $\ep$ small enough, the solution
$(X^\ep_t, Z^\ep_t)$ of $\widetilde {(M+D)_\ep}$ is also a solution of
the BSDE $\widetilde {(M+D)}$ on a small time interval $[T^\ep_1;T]$.
This result relies on Proposition \ref{zbounded}, which gives the very
strong condition that $Z$ is bounded a.s.

Firstly, let us give a classical link between BSDEs and PDEs.

\begin{sprop}
Let $\ep \in ]0;1[$; we use the notations introduced in paragraph
\ref{par4.2} (recall in particular that the function $\gamma$ depends
on $\ep$).

\noi
Consider for
$u_\ep=(u_\ep^1, \ldots, u_\ep^n) :
          [0;T] \times \R^d \rightarrow \R^n $
the following system of quasilinear parabolic partial differential
equations
\begin{equation}
\left\{
\begin{array}{ccl}
    \frac{\partial u_\ep}{\partial t}(t,x)
 & = & {\cal L} u_\ep (t,x)
    +\gamma(x, u_\ep(t,x), (\nabla_x u_\ep \si)(t,x)) \\
  u_\ep (0,x)  &=& F(x)
\end{array}
\right.
\label{eqchaleur}
\end{equation}
where $\nabla_x u_\ep$ is the $n \times d$ matrix whose rows
are $(\nabla_x u_\ep^i)_{i=1, \ldots, n}$, the partial derivatives of
the components of $u_\ep$ with respect to space; moreover,
\begin{equation}
{\cal L} u_\ep = (L u_\ep^1, \ldots, L u_\ep^n)
\label{difopL1}
\end{equation}
is a vector in $\R^n$, and
\begin{equation}
 L = \half \sum_{i,j=1}^d (\si {}^t \si)_{i,j}(t,x)
   \frac{\partial^2}{\partial x_i \partial x_j} + \sum_{i=1}^d b_i(t,x)
   \frac{\partial}{\partial x_i}
\label{difopL2}
\end{equation}
is the infinitesimal generator of the diffusion \eqref{sde}.

Then
\begin{item}{(i)}
This equation has a unique solution $u_\ep$
in $C^{1,2}([0;T] \times \R^d,\R^n)$ (i.e. $u_\ep$ has continuous first
derivative with respect to time and continuous second derivatives with
respect to space);
\end{item}

\begin{item}{(ii)}
A.s., $X_t^\ep = u_\ep (T-t,B_t^y) \ \ \forall t$;
\end{item}

\begin{item}{(iii)}
$\forall t$, a.s.,
$Z_t^\ep = \nabla_x u_\ep (T-t,B_t^y) \sigma(B_t^y)$.
\end{item}
\end{sprop}

\begin{demo}
The assertion $(i)$ is a classical result; for a probabilistic proof,
see Theorem 3.2 of \cite{pp92}, and for a proof in dimension 1 (i.e.
$n=1$) with less regular functions $\gamma$, see \cite{bl97}.

Then $(ii)$ and $(iii)$ follow easily~:
both It\^o's formula and equation \eqref{eqchaleur} lead to
$$
 \left\{
  \begin{array}{ccl}
d u_\ep(T-t,B_t^y)
& = & -\gamma(B_t^y, u_\ep, (\nabla_x u_\ep \si)(B_t^y)) dt
                 + (\nabla_x u_\ep \si) d W_t   \\
u_\ep(0,B_T^y)
& = & F(B_T^y) = U.
  \end{array}
 \right.
$$

Thus
$$\left( u_\ep(T-t,B_t^y), \nabla_x u_\ep (T-t, B_t^y) \sigma(B_t^y)
   \right) _t$$
is a solution of BSDE $\widetilde {(M+D)_\ep}$. But $(X^\ep,Z^\ep)$ is
also a solution and this BSDE has Lipschitz coefficients, so it has a
unique solution. Then $(ii)$ and $(iii)$ follow.
 \end{demo}

{\bf Remark. }Since $\nabla_x u_\ep$ is continuous, the trajectories
$t \mapsto Z_t^\ep$ can be taken continuous.

\bigskip
Now we give the main result of this paragraph.

\begin{sprop}
There is an $\ep \in ]0;1[$  and a $T^\ep_1 \in [0;T[$ (depending on
$\ep$) such that a.s. for $t \in [T^\ep_1;T]$, we have
$\Vert Z_t^\ep \Vert \le \frac{1}{\ep}$. This means, with the
notations of paragraph \ref{par4.2}, that $h_\ep(Z_t^\ep)=0$ so
$\overline Z_t^\ep = Z_t^\ep$ and $(X_t^\ep,Z_t^\ep)$ is a solution
of BSDE $\widetilde {(M+D)}$ on the time interval $[T^\ep_1;T]$.
\label{zbounded}
\end{sprop}

\begin{demo}
We only deal with a drift
$f$ depending both on $x$ and $z$, the $z$-independent case being
similar (and easier).

We let $\al_\ep$ be the global Lipschitz constant of the function
$\gamma$ (independent of $b$):
\begin{eqnarray}
  & & \forall (x,z), (x',z') \in \R^n \times \R^{n d_W}, \nonumber  \\
  & & \vert \gamma(b,x,z) - \gamma(b,x',z') \vert \le
 \al_\ep \left( \vert x - x' \vert + \Vert z - z'
                                       \Vert \right).
\label{lipphi}
\end{eqnarray}
(in general, $\al_\ep$ will tend to $\infty$ as $\ep$ decreases
to zero).

We have to show that
$\Vert Z_t^\ep \Vert = \left\Vert \nabla_x u_\ep
    (T-t,B_t^y) \sigma(B_t^y) \right\Vert \le \frac{1}{\ep}$ a.s.,
so it suffices to check that $u_\ep$ is
$\frac{1}{\ep \Vert \sigma \Vert}_\infty$-Lipschitz
with respect to the space variable (${\Vert \si \Vert}_\infty$ denoting
the supremum of $\Vert \si(b) \Vert$ for $b \in \R^d$). This work will
be achieved in several steps.

{\it First Step : }
We consider two solutions
$(Y,Z)$ and $(\hY, \hZ)$ of $\widetilde {(M+D)_\ep}$ corresponding to two
terminal values $Y_T$ and $\hY_T$.
We set $\dY_t = Y_t - \hY_t$ and $\dZ_t = Z_t - \hZ_t$.

Then
$$
  \dY_t + \int_t^T \dZ_s d W_s =
  \dY_T + \int_t^T \left(\gamma (B_s^y,Y_s,Z_s)
   - \gamma (B_s^y, \hY_s,\hZ_s) \right) ds.
$$

Independence gives :
$$
  \esp \vert \dY_t \vert ^2 + \esp \left(
   \int_t^T \Vert \dZ_s \Vert ^2 d s \right) =
   \esp \left( \left\vert \dY_T + \int_t^T
  (\gamma (B_s^y, Y_s,Z_s) - \gamma (B_s^y,\hY_s,\hZ_s))
   d s \right\vert ^2 \right).
$$

Using H\" older's inequality, we obtain :
$$
\begin{array}{rcll}
  \esp \vert \dY_t \vert ^2 + \esp \left(
     \int_t^T \Vert \dZ_s \Vert ^2 ds \right)   \le
& 2 & \esp \vert \dY_T \vert ^2 &                               \\
& + & 2(T-t) \esp
       \Bigg(
&  \int_t^T
      \vert \gamma (B_s^y,Y_s,Z_s)       \\
&   & &    - \gamma (B_s^y,\hY_s,\hZ_s) \vert  ^2 ds \Bigg).
\end{array}
$$
Hence using \eqref{lipphi} we have
$$
\begin{array}{rcll}
  \esp \vert \dY_t \vert ^2 + \esp \left(
   \int_t^T \Vert \dZ_s \Vert ^2 d s \right) \le
& 2 & \esp \vert \dY_T \vert ^2 &                           \\
& + & 4 \alpha_\ep^2 (T-t)
   \Bigg(
 & \esp \left( \int_t^T \vert \dY_s \vert ^2 d s \right)   \\
&   & & +\esp \left( \int_t^T \Vert \dZ_s \Vert ^2 d s \right) \Bigg).
\end{array}
$$

{\it Second Step : }
We want to remove the terms
$\esp \left( \int_t^T \Vert \dZ_s \Vert ^2 d s \right)$, to get
an equation with $Y$ only; this leads to make the
following assumption (which will be true in the end of the proof) :
$t \in [T^\ep_1;T]$ with $T^\ep_1 \in [0;T[$ such that
$4 \alpha_\ep^2 (T-T^\ep_1) \le 1$.

Then, if $t \in [T^\ep_1;T]$, we get
$$
      \esp \vert \dY_t \vert ^2
 \le
      2 \esp \vert \dY_T \vert ^2
      + \esp \left( \int_t^T \vert \dY_s \vert ^2 ds \right).
$$
And Gronwall's lemma gives
\begin{equation}
\esp \vert \dY_t \vert ^2
           \le 2 e^{T-t} \esp \vert \dY_T \vert ^2.
\label{gronwall}
\end{equation}

{\it Third Step : }
Let us choose $Y_T=F(B^{t,x}_T)$ and $\hY_T=F(B^{t,\hat x}_T)$
(where $(B_.^{t,x})$ denotes the diffusion starting at $x$ at time $t$);
then
\begin{equation}
\dY_t= u_\ep(T-t,x) - u_\ep(T-t, \hat x)
\label{y=u}
\end {equation}
and
\begin{equation}
\esp \vert \dY_T \vert ^2 \le L_F^2 \esp \left(
     \vert B_T^{t,x} - B_T^{t,\hat x} \vert^2 \right)
     \le L_F^2 C_{\si, b} \vert x - \hat x \vert ^2,
\label{majoryT}
\end{equation}
with $L_F$ the Lipschitz constant of $F$, and the last inequality is
a well-known result of $L^2$-continuity with respect to initial
conditions. Hence from \eqref{gronwall}
$$ \vert u_\ep(T-t,x) - u_\ep(T-t, \hat x) \vert \le
   L_F \sqrt{2C_{\si, b} e^{T-t}} \vert x - \hat x \vert$$
and
$$ \left\Vert \nabla_x u_\ep (T-t,x) \right\Vert
      \le    L_F \sqrt{2C_{\si, b} e^T}.$$
(In fact, the right hand side is given up to a positive multiplicative
constant, due to the different norms used; it doesn't matter in the
sequel).

{\it Conclusion : }
Let $\ep \in ]0;1[$ such that
$\frac{1}{\ep} \ge \Vert \sigma \Vert_\infty L_F \sqrt{2C_{\si, b} e^T}$.
Then the preceding
proof shows that if $T^\ep_1 \in [0;T[$ is such that
$4 \al_\ep^2 (T-T^\ep_1) \le 1$, then for $t \in [T^\ep_1;T]$
we have
\begin{equation}
 a.s. \ \ \Vert Z_t^\ep \Vert =
       \left\Vert \nabla_x u_\ep
       (T-t,B_t^y) \sigma(B_t^y) \right\Vert  \le  \frac{1}{\ep}.
\label{zinfep}
\end{equation}
In fact, using the continuity in $t$ of $Z_t^\ep$, this inequality holds
for any $t$, a.s.; that finishes the proof of the proposition.
\end{demo}

\subsection{A solution of $(M+D)$ on a small time interval}
\label{par4.4}

The framework is the same as the one introduced in Subsections
\ref{par4.2} and \ref{par4.3}.
The aim of this section is to prove that the above solution
$(X_t,Z_t)_t$ (we omit in this paragraph the superscript $\ep$ for
notational convenience) of $\widetilde {(M+D)}$ is, under an
additional condition on the drift $f$, a solution of
$(M+D)$ (in fact we prove that $(X_t)$ remains in $\omb$, the compact
where the terminal value $U$ lies).
We recall that $\omb= \{ \chi \le c\}$ is the sublevel set of a smooth
convex (for the connection $\Gamma$) function $\chi$, defined on an
open set $O$ relatively compact in $\R^n$, and moreover that $\omb$
is relatively compact in $O$. We make the following hypothesis

$$(H_s) \ \ f \hbox{ is pointing strictly outward on the boundary }
\partial \omb \ \hbox{ of } \omb.$$
It means that
\begin{equation}
 \forall (b,x,z) \ : \ x \in \partial \omb, \ \
   \underset{b,x,z}{\inf}(D \chi(x) | f(b,x,z))_r \ge \zeta >0,
\label{Hs}
\end{equation}
where $(\cdot | \cdot)_r$ denotes the Riemannian metric tensor (if
$f^\bot$ is the component of $f$ orthogonal to
$\partial \omb = \{ \chi = c\}$, it is equivalent to require that
$\underset{b,x,z}{\inf} (D \chi(x) | f^\bot (b,x,z))_r$ be bounded
below by a positive constant, since
$\chi$ is constant (equal to $c$) on $\partial \omb$).

{\bf Remark :} This condition arises naturally in the deterministic
version of equation $(M+D)$
$$
\left\{
      \begin{array}{l}
          d x_t =  f(x_t) d t            \\
	    x_T = u                      \\
      \end{array}
    \right.   $$
where $u$ is deterministic and so $Z_t=0$ for any $t$ (in fact, as we
shall see in Subsection \ref{par4.6}, the natural condition is to
require that the infimum in \eqref{Hs} be only nonnegative).

\begin{sprop}
\label{prop4.1}
Suppose that $\chi$ is strictly convex on $O$ (this means that
$\Hess \chi$ is positive definite). Then, under the
assumption $(H_s)$, the process $(X_t)_{t \in [T^\ep_1;T]}$ remains in
$\omb$, i.e. $(X_t,Z_t)$ is a solution of $(M+D)$ on the time interval
$[T^\ep_1;T]$.
\end{sprop}

\begin{demo}
Once again, the goal is to construct a nonnegative submartingale,
null at time $T$, which vanishes if and only if the process $X$ is in
$\omb$. The proof will be split into three steps.

{\it First Step : }Framework of the proof.

Suppose that $c>0$ and $\chi$ reaches its infimum at $p \in \omega$
with $\chi(p)=0$. Then consider the following mapping, defined on a
normal open neighbourhood of $\omb$ centered at $p$ (so it is in
particular a neighbourhood of $0$ in $\R^n$)
$$ y \mapsto \frac{\sqrt{c_2} y}
                  {\sqrt{c_2 \Vert y \Vert^2 + c-\chi(y)} };$$
for $c_2>0$ small enough, it is a diffeomorphism from an open set
$O_1$ (relatively compact in $O$ and containing $\omb$) onto an open
neighbourhood $N$ of $\overline {B(0,1)}$, such that $\omb$ is sent
onto $\overline {B(0,1)}$ (in fact, it is sufficient to take
$c_2 \le \half \lambda_\chi$, where $\lambda_\chi$ denotes the
(positive) infimum on $O_1$ of the eingenvalues of $\Hess \chi$).
Using this diffeomorphism, we can work in a local chart $O$ (take
$O:=N$) such that $\omb = \overline {B(0,1)}$ and $\chi(0)=0$.

{\it Second Step : }Construction of a "nearly convex" function $H$.

Choose $\rho>1$ such that $\overline {B(0,\rho)} \subset O$; then we
take for the mapping $\phi$ (see the beginning of Subsection \ref{par4.2})
a smooth function $\R^n \rightarrow [0;1]$ equal to $1$
on $B(0,1)$ and to $0$ outside $B(0,\rho)$, which has moreover
spherical symmetry.
This gives $\tilde \Gamma=\phi \Gamma$ and $\tilde f=\phi f$.
Then define on $\R^n$
$$k(x)=\phi(x) \chi(x)+(1-\phi(x)) \alpha(x)$$
where
$\alpha(x)=a \vert x \vert$ and $a>0$ is chosen so that $k(x) \le c$
if and only if $x \in \overline {B(0,1)}$ (take for instance
$a > \sup_{B(0,\rho)} \chi$);
$k$ is clearly a nonnegative smooth mapping.
We would like to have a mapping which vanishes on $\omb$, so we
let $H=h \circ k$ where $h : \R_+ \rightarrow \R_+$ is a smooth convex
(so nondecreasing) function, vanishing on the interval $[0;c]$ (only)
and growing linearly at infinity.

The mapping $H$ so defined is convex for the connection
$\tilde \Gamma$ on $B(0,1)$ and outside $B(0,\rho)$. Indeed, on
$B(0,1)$, $H=0$ and outside $B(0,\rho)$, $H=h \circ \alpha$
which is convex for the flat connection ($=\tilde \Gamma$).

{\it Third Step : }
We work on the time interval $[T^\ep_1;T]$. The aim of this step is to
show that, for $\la$ large enough, the process
$(e^{\lambda t} H(X_t))_t$ is a real submartingale. Reasoning as in the
uniqueness part, we apply It\^o's formula~:
\begin{eqnarray*}
e^{\lambda t} H (X_t) - e^{\lambda T_1^\ep} H (X_{T_1^\ep})
& = &\int_{T_1^\ep}^t e^{\lambda s} D H (X_s) (Z_s d W_s)    \\
&   &+ \half \int_{T_1^\ep}^t e^{\lambda s}
    \left(  \sum_{i=1}^{d_W}  {}^t [{}^t Z_s]^i
                       \widetilde {\Hess} H (X_s)
                   [{}^t Z_s]^i \right) d s                  \\
&   & + \int_{T_1^\ep}^t e^{\lambda s} D H (X_s)
                          \cdot \tilde f (B_s^y,X_s,Z_s) d s \\
&   & + \int_{T_1^\ep}^t \lambda e^{\lambda s} H (X_s) d s.
\end{eqnarray*}

The stochastic integral is a martingale because $D H$ is
bounded; it remains to prove that the bounded variation term is an
increasing process, i.e. show the nonnegativity of the sum
$$
\half
  \sum_{i=1}^{d_W} {}^t [{}^t z]^i \widetilde{\Hess} H (x) [{}^t z]^i
    +  D H(x) \tilde f(b,x,z) + \lambda H(x).
$$

But $\widetilde{\Hess} H (x) \ge h'(k(x)) \Hess k (x)$ since $k$
is nondecreasing (see for instance (4.36) p 42 in \cite{scm}); then
it suffices to prove the nonnegativity of

\begin{equation}
h'(k(x)) \left(
\half
  \sum_{i=1}^{d_W} {}^t [{}^t z]^i \widetilde{\Hess} k (x) [{}^t z]^i
    +  D k(x) \tilde f(b,x,z) \right) + \lambda H(x).
\label{pos3}
\end{equation}

Remark that $h'(k(x)) \ge 0$ and that, because of the boundedness of
the process $Z_t$ (according to \eqref{zinfep}), it is sufficient to
consider $z$ such that $\Vert z \Vert \le \frac{1}{\ep}$.

Let $A_1$, $A_2$ and $A_3$ denote the three terms in this order
in the sum \eqref{pos3}.

\noindent $\bullet$
if $x \in \overline {B(0,1)}$, $A_1$ and $A_2$
vanish because $H=0$, so the sum is nonnegative;

\noindent $\bullet$
if $x \in {}^cB(0,\rho)$, $A_1 \ge 0$ because $k=\alpha$
is convex for the flat connection and $A_2=0$ because $\tilde f=0$;
then the sum is nonnegative.

\noindent $\bullet$
if $x \in B(0,\rho) \setminus \overline {B(0,1)}$, we
consider two situations :
firstly when $x$ belongs to a neighbourhood of the sphere $S(0,1)$
(i.e. when $\vert x \vert \in ]1;1+\eta[$ where $\eta>0$ is to be
determined). Using continuity, $\widetilde {\Hess} k (x)$ is
nonnegative (in the sense of matrices) on $]1;1+\eta_1[$
(with $\eta_1>0$ sufficiently small) because if $\vert x \vert =1$,
$\widetilde {\Hess} k (x)=\widetilde {\Hess} \chi (x)$ and
$\chi$ is strictly convex on $B(0,1)$ for the connection
$\tilde \Gamma$. This gives $A_1 \ge 0$. \\
For $A_2$, using hypothesis $(H_s)$, we get that for
$x \in \partial \omb$ and $\Vert z \Vert \le \frac{1}{\ep}$,
$$ Dk(x) \cdot \tilde f(b,x,z) = D \chi(x) \cdot \tilde f (b,x,z)
                   = D \chi(x) \cdot f (b,x,z) \ge \zeta >0.$$
So by uniform continuity and hypothesis \eqref{lip} and
\eqref{upperboundf}, there is an $\eta_2>0$ such that
$$\forall b, \
  \forall y : 1< \vert y \vert < 1+\eta_2, \
  \forall z : \Vert z \Vert \le \frac{1}{\ep}, \
 D k(y) \cdot \tilde f (b,y,z) \ge 0.$$
In particular, $A_2 \ge0$ for $1<\vert x \vert < 1+\eta_2$.
Let $\eta = \min(\eta_1, \eta_2, \rho-1)$. Then $A_1+A_2 \ge 0$ if
$1<\vert x \vert < 1+ \eta \le \rho$,
and the sum \eqref{pos3} is nonnegative.

On the other hand, if $\vert x \vert \in ]1+\eta;\rho[$ (i.e.
if $\vert x \vert$ is "far" from $1$), then $A_1+A_2$ is
bounded above since $H$ is smooth, and $\tilde f (b,x,z)$ is
bounded
since $z$ is bounded. But we have constructed $H$ so that
$H(x) \ge \theta >0$ if $\vert x \vert >1+\eta$.
Then we can choose $\lambda>0$ such that $A_1+A_2+A_3$ is
nonnegative.

{\it Conclusion : }
The end of the proof now goes on by classical arguments~:
the process $(e^{\lambda t} H(X_t))_{t \in [T_1^\ep;T]}$ is a
nonnegative submartingale for the $\lambda$ chosen above, null at
time $T$. Then $\forall t \in [T_1^\ep;T], \ H(X_t)=0$; and the
mapping $H$ has precisely been chosen so that this implies
$X_t \in \omb$.
\end{demo}

{\bf Remark :} In our local coordinates, let $D_r \chi (x)$ denote
the radial derivative of $\chi$ at $x$ and $f^r$ the radial component
of $f$. Then, since $\chi = c$ on the sphere $S(0,1)= \partial \omb$,
we have
$$ \forall x \in \partial \omb, \ D \chi (x) \cdot f(b,x,z)
   = D_r \chi (x) \cdot f^r(b,x,z).$$
But if $\Hess \chi$ is positive definite, $D_r \chi(x)$ is a positive
real number and in our coordinates, $(H_s)$ is equivalent to require
that $f^r(b,x,z)$ be bounded below by a positive constant, independent
of $x \in \partial \omb$ and of $b,z$.
\medskip

Now we come back to the end of the proof of the existence.

\subsection{The solution on the whole interval $[0;T]$}
\label{par4.5}

According to Subsections \ref{par4.3} and \ref{par4.4},
$(X_t^\ep,Z_t^\ep)$ is a solution of BSDE $(M+D)$ on the time interval
$[T_1^\ep;T]$; moreover, $\ep$ and $T_1^\ep$ must verify (see the
conclusion in the proof of Proposition \ref{zbounded})
$$
\left\{
\begin{array}{rcl}
\frac{1}{\ep}           & \ge & \Vert \sigma \Vert_\infty L_F
                           \sqrt{2 C_{\si, b} e^T}  \\
4 \al_\ep^2 (T-T_1^\ep) & \le & 1
\end{array}
\right.
$$
where $\alpha_\ep$ tends {\it a priori} to $\infty$ when
$\ep$ goes to zero. The aim of this subsection is to show that it is
a solution on the time interval $[0;T]$.

In general, if we take back the proof of Proposition \ref{zbounded},
with now $T_1^\ep$ and $u_\ep(T-T_1^\ep, B_{T_1^\ep}^y)$
as terminal time and variable, we get a solution of $(M+D)$ on a time
interval $[T_2^\ep;T_1^\ep]$;
but it is easy to show that $\ep$ and $T_2^\ep$ must verify now
$$
\left\{
\begin{array}{rcl}
\frac{1}{\ep} & \ge & (\Vert \sigma \Vert_\infty  L_F
               \sqrt{2 C_{\si, b} e^T}) \sqrt{2e^T} \\
4 \al_\ep^2 (T_1^\ep-T_2^\ep) & \le & 1;
\end{array}
\right.
$$
as a consequence, the length of the interval $[T_2^\ep;T_1^\ep]$ may
be less than the one of $[T_1^\ep;T]$ and repeating this method
inductively could lead to a solution on an interval $]T_0;T[$ with
$T_0>0$ only (this means that the solution explodes).

\bigskip
In fact, the existence of $\Psi$ prevents the solution from exploding
and allows to build inductively a solution of $(M+D)$ on $[T_1^\ep;T]$,
$[T-2(T-T_1^\ep);T]$, $[T-3(T-T_1^\ep);T]$, ... and so to get at the
end a solution on $[0;T]$.
\begin{sprop}
Suppose that assumption $(H_s)$ holds and that $\chi$ is strictly convex
on $O$.
If $\ep>0$ is small enough, then $(X_t^\ep,Z_t^\ep)$, the solution of
BSDE $\widetilde{(M+D)}_\ep$, is also a solution
of BSDE $(M+D)$ on the whole time interval $[0;T]$.
\end{sprop}

\begin{demo}
Consider $T_1^\ep$ and $u_\ep(T-T_1^\ep, B_{T_1^\ep}^y)$
as terminal time and variable; applying the two preceding sections, we
get a solution of $(M+D)$ (with $X_t^\ep \in \omb$) on a time
interval $[T_2^\ep;T_1^\ep]$ where $\ep$ and $T_2^\ep$ verify
$$
\left\{
\begin{array}{rcl}
\frac{1}{\ep} & \ge & \Vert \sigma \Vert_\infty   L_{u_\ep}         \\
4 \al_\ep^2 (T_1^\ep-T_2^\ep) & \le & 1
\end{array}
\right.
$$
($L_{u_\ep}$ is the Lipschitz constant of $u_\ep$ for the
space variable, uniformly in $t$ on $[T_2^\ep;T]$).

As in Subsection \ref{par4.3},
for $t \in [T_2^\ep;T[$ we let $(Y_s,Z_s)_{s \in [t;T]}$ and
$(\hat Y_s,\hat Z_s)_{s \in [t;T]}$ be two solutions of BSDE
$\widetilde{(M+D)_\ep}$ such that
$Y_T=F(B_T^{t,x})$ and $\hat Y_T=F(B_T^{t,\hat x})$ (by which
we denote diffusions starting at $x$ or $\hat x$ at time $t$).
According to Subsection
\ref{par4.4}, these two processes remain in $\omb$, so we can make
use of the function $\Psi$. Then the same inequalities as in
\eqref{major2} give

\begin{eqnarray*}
  \esp \left( \Psi(Y_t, \hat Y_t) \right)
& \le &
  C_1 \esp \left( e^{\eta \int_0^T (\Vert Z_u \Vert_r^2 +
  \Vert \hat Z_u \Vert_r^2) du} \right)^\frac{1}{4}
  \esp \left( \delta^{4p}(Y_T, \hat Y_T) \right)^\frac{1}{4}        \\
& \le &
  C \esp \left( \delta^{4p}(Y_T, \hat Y_T) \right)^\frac{1}{4}      \\
\end{eqnarray*}
(the second inequality is obtained by letting $\eta=0$ in the
$z$-independent case; in the other case, by letting $\eta=\al$ as in
Corollary \ref{unifinteg} and $p=2$).

The equivalence of the Riemannian and Euclidean distances on
$\omb \times \omb$ and $\Psi \approx \delta^p$, together with
\eqref{y=u} and \eqref{majoryT} give
$$ \vert u_\ep(T-t,x)-u_\ep(T-t,\hat x) \vert
              \le C_0 L_F \vert x- \hat x \vert,$$
where $C_0$ depends only on the diffusion $(B_t^y)_t$, $T$, $\Psi$,
$\omb$ and the drift $f$ (then $L_{u_\ep} = L_F C_0$).

Consequently, if we take $\ep$ such that
\begin{equation}
\frac{1}{\ep} \ge \Vert \sigma \Vert_\infty L_F C_0,
\label{condepsilon}
\end{equation}
$T_1^\ep$ s.t. $4 \al_\ep^2 (T-T_1^\ep) \le 1$ and $T_2^\ep$ s.t.
$4 \al_\ep^2 (T_1^\ep-T_2^\ep) \le 1$, we get a solution on an interval
$[T_2^\ep-T_1^\ep]$ with the same length as $[T_1^\ep;T]$.
Repeating the same method inductively with the same $\ep$ at each step,
we obtain that $(X_t^\ep,Z_t^\ep)$ is a solution of $(M+D)$ on
$[T_1^\ep;T]$, $[T-2(T-T_1^\ep);T]$, $[T-3(T-T_1^\ep);T]$, ... and so
on the whole interval $[0;T]$.
\end{demo}

{\bf Remark :} As $Z_t^\ep = \nabla_x u_\ep(T-t,B_t^y) \sigma(B_t^y)$,
it is a
straightforward consequence of the proof that a.s., for all $t$,
$\Vert Z_t^\ep \Vert$ is bounded above by $\frac{1}{\ep}$.

If we sum up the results obtained, we have the

\begin{sprop}
\label{existence1}
We consider BSDE $(M+D)$
with a terminal value $U$ in $\omb= \{ \chi \le c \}$. Suppose that $f$
is a $C^3$ function which verifies conditions \eqref{lip},
\eqref{upperboundf} and $(H_s)$.
If moreover $\chi$ is strictly convex (i.e. $\Hess \chi$ is positive
definite), then
\begin{item}
(i) If $f$ does not depend on $z$, the BSDE has a (unique) solution
such that $X$ remains in $\omb$.
\end{item}
\begin{item}
(ii) If $M$ is a Cartan-Hadamard manifold and the Levi-Civita connection
is used, then the BSDE has a (unique) solution such that $X$ remains in
$\omb$ too.
\end{item}
\end{sprop}

The last paragraph is devoted to generalize Proposition \ref{existence1}
to drifts $f$ which are less regular and satisfy a weaker hypothesis
than $(H_s)$.

\subsection{The solution for general $f$}
\label{par4.6}

Let $f$ be a function verifying \eqref{lip}, \eqref{upperboundf} and
the hypothesis
$$(H) \ \ f \hbox{ is pointing outward on the boundary of } \omb$$
introduced in Subsection \ref{par1.4}.
This means that
$$
 \forall (b,x,z) \ : \ x \in \partial \omb, \ \
   (D \chi(x) | f(b,x,z))_r \ge 0;
$$
equivalently, if $f^\bot$ is the component of $f$ orthogonal to
$\partial \omb = \{ \chi = c\}$, we may require
$(D \chi(x) | f^\bot (b,x,z))_r$ to be nonnegative.
This condition is obviously weaker than $(H_s)$.

Firstly, we will derive from equation $(M+D)$ new BSDEs, each one having
a unique $\omb$-valued solution. Then we will show that these solutions
converge to the solution of $(M+D)$. The function $\chi$ will supposed to
be strictly convex and calculus will be achieved for a drift $f$ depending
both on $x$ and $z$ (the case when $f$ does not depend on $z$ being
simpler).

\bigskip
Let us consider the
local chart $O$ defined in the first step of the
proof of Proposition \ref{prop4.1}. Remember that
$\omb=\overline{B(0,1)}$ in these local coordinates and remark that
hypothesis $(H)$ means that the radial component $f^r(b,x,z)$ of
$f(b,x,z)$
is nonnegative for $x \in \partial \omb$ (see the remark at the end of
Subsection \ref{par4.4}).
Extend the mapping $f$ to $\R^d \times \R^n \times \R^{n d_W}$ by putting
$f(b,x,z)=0$ if $x \notin O$ and define (on
$\R^d \times \R^n \times \R^{n d_W}$) the
convolution product for $l \in \N^*$
$f_l = f * \rho_l$
where $\rho_l(b,x,z) = l \rho (l \Vert (b,x,z) \Vert)$ and
$\rho : \R_+ \rightarrow \R_+$ is a bump
function (i.e. a smooth function with $\rho'(0)=0$, $\rho=0$ outside
$[0;1]$ and $\int_{\R_+} \rho(u) du=1$).
Besides, let us define a function $g_l$ by
\begin{equation}
g_l(b,x,z) = f_l(b,x,z)  + \frac{A}{l} x,
\label{defgl}
\end{equation}
where $A$ is a positive constant which will be chosen below; we introduce
on $O$ the following BSDE for $F \in C_c^\infty(\R^d)$ and $U \in \omb$
$$(M+D)_l
\left\{
      \begin{array}{l}
          d X_t = Z_t d W_t + \left( - \half \Gamma_{jk}(X_t)
            ([Z_t]^k \vert [Z_t]^j) + g_l(B_t^y,X_t,Z_t) \right) d t \\
	  X_T =  F(B_T^y) = U.                                      \\
      \end{array}
    \right.   $$

\begin{slemma}
Let $O_1$ be an open set such that
$\omb \subset O_1 \subset \overline{O_1} \subset O$. Then for $l \ge l_0$
the functions $g_l$ are smooth on $\R^d \times O_1 \times \R^{n d_W}$ and
verify a
Lipschitz condition like \eqref{lip2} with the same Lipschitz constant
$L''$. Moreover, they also verify a boundedness condition like
\eqref{upperboundf}.
\label{uniflip_gl}
\end{slemma}

\begin{demo}
For $x \in O_1$, we write
$$ f_l(b,x,z)
    = \int_{\R^d \times \R^n \times \R^{n d_W}}
                f((b,x,z)-(\beta,y,w)) \rho_l(\beta,y,w) d (\beta,y,w).$$
Thus, as soon as $dist(O_1,O) > 1/l_0 \ge 1/l$, $\rho_l(\beta,y,w)=0$ if
$\vert y \vert \ge dist(O_1,O)$; so the integrand vanishes if
$x-y \notin O$ and we can use the Lipschitz property \eqref{lip2} of $f$
on $O$. Then the properties of convolution give the result for $f_l$. As
it also holds obviously for the functions $x \mapsto A/l \ x$, we have the
result for $g_l$.
The second assertion is an easy consequence of conditions \eqref{lip}
and \eqref{upperboundf} for the drift $f$.
\end{demo}

In the sequel, we will consider the sequence $(g_l)_l$ only for
$l \ge l_0$.

\begin{sprop}
For every $l$, the BSDE $(M+D)_l$ has a (unique) $\omb$-valued solution
$(X^l,Z^l)$; moreover, there is an $\ep>0$, independent of $l$, such that
a.s., $\Vert Z_t^l \Vert \le \frac{1}{\ep}$ for any $t$.
\end{sprop}

\begin{demo}
Using Lemma \ref{uniflip_gl}, we apply Subsections \ref{par4.2} and
\ref{par4.3} to $(M+D)_l$. We get an $\ep_l$ and a $T_1^{\ep_l}$ for every
$l$; but the Lipschitz constant of $\gamma$ is independent of $l$ (since
the $g_l$ have the same one), so the proof of Proposition \ref{zbounded}
shows that we can choose $\ep_l=\ep$ and $T_1^{\ep_l}= T_1^\ep$,
independently of $l$.

In order to apply Subsection \ref{par4.4}, we need to prove that $g_l$
verifies condition $(H_s)$; in fact, we have seen (see in particular the
remark at the end of Subsection \ref{par4.4}) that it suffices to show
that
\begin{equation}
 \forall (b,x,z) \ : \ x \in \partial \omb, \ \Vert z \Vert \le
    \frac{1}{\ep}, \ \
   \underset{b,x,z}{\inf} g_l^r(b,x,z) \ge \zeta >0,
\label{glrpos}
\end{equation}
It is easy to see with the properties of convolution
that, for a constant $\hat C$ depending only on $\ep$ and $f$,
$$
\forall l, \ \forall b, \ \forall x \in \partial \omb, \ \forall z :
\Vert z \Vert \le
    \frac{1}{\ep}, \ \vert f_l(b,x,z)-f(b,x,z) \vert \le \frac{\hat C}{l}.
$$
This and the nonnegativity of $f^r(b,x,z)$ for $x \in \partial \omb$ give
$$
\forall l, \ \forall b, \ \forall x \in \partial \omb, \ \forall z :
\Vert z \Vert \le
    \frac{1}{\ep}, \ f_l^r(b,x,z) \ge - \frac{\hat C}{l}.
$$
Now if we take in \eqref{defgl} $A=\hat C + 1$, obviously \eqref{glrpos}
holds with $\zeta = 1/l$ and Subsection \ref{par4.4} can be applied.

Then the results of Subsection \ref{par4.5} hold (in particular
Proposition \ref{existence1}) with the same $\ep$ for every $l$ (this
comes from \eqref{condepsilon}, remarking that the constant $C_0$ in
this inequality is independent of $l$).
The proof is completed.
\end{demo}

Now we prove that, as expected, the limit of these solutions when
$l \rightarrow \infty$ solves equation $(M+D)$.

\begin{sprop}
For $n \in \N$ and $s \in [0;T]$, we note for simplicity
$V_s^n := (B_s^y,X^n_s,Z^n_s)$,
When $l$ tends to infinity, the above solution $(X^l, Z^l)$ converges
(for the usual $L^2$ norms for $X$ and $Z$) to a pair $(X,Z)$ which
solves equation $(M+D)$.
\end{sprop}

\begin{demo}
For $l$ and $m$ in $\N$, we note
$\dX_t=X_t^l-X_t^m$, $\dZ_t=Z_t^l-Z_t^m$ and

$$ A_{l,m} = \esp \int_0^T \vert \dX_s \vert
   \left( \vert f_m(V_s^l) - f_l (V_s^l) \vert
   + A \left\vert \frac{1}{l} X_s^l
                - \frac{1}{m} X_s^m \right\vert \right) ds.$$

Applying It\^o's formula (in $\R^n$) to $\vert \dX_t \vert^2$ between
$t$ and $T$, we get (note that a dot stands for the inner product in
$\R^n$)

\begin{eqnarray}
- \vert \dX_t \vert^2
& = & \int_t^T 2 \dX_s \cdot (\dZ_s dW_s) + \int_t^T 2 \dX_s \cdot
      \left( g_l(V_s^l)- g_m(V_s^m) \right)
                                                        ds  \nonumber\\
&   & - \half \int_t^T 2 \dX_s \cdot \left( \Gamma_{jk}(X^l)
      ([Z^l_s]^k \vert [Z^l_s]^j) - \Gamma_{jk}(X^m)
      ([Z^m_s]^k \vert [Z^m_s]^j) \right) ds               \nonumber \\
&   & + \int_t^T \Vert \dZ_s \Vert^2 ds.
\label{itoRn}
\end{eqnarray}

Using the uniform boundedness of the $(Z^l_t)$ (proved in the
preceding lemma), we can bound above the integral involving the
Christoffel symbols by $C \int_t^T \vert \dX_s \vert ^2 ds$; besides,
for the second term on the right, we write
\begin{eqnarray*}
\left\vert \int_t^T  \dX_s \cdot
 (g_l(V_s^l)- g_m(V_s^m)) ds \right\vert
& \le & A_{l,m} + \int_t^T \vert \dX_s \vert \cdot
        \vert f_m(V_s^l)- f_m(V_s^m) \vert ds       \\
& \le & C_1 \left( A_{l,m} + \int_t^T \vert \dX_s \vert
         (\vert \dX_s \vert + \vert \dZ_s \vert) ds \right)     \\
& \le & C \left( A_{l,m} + \int_t^T \vert \dX_s \vert^2 ds \right)
           + \half \int_t^T \Vert \dZ_s \Vert^2 ds
\end{eqnarray*}
where $C$ is independent of $l$ and $m$. Then, we obtain by taking
the expectation in \eqref{itoRn}

\begin{equation}
\esp \vert \dX_t \vert^2 + \half \esp \int_t^T \Vert \dZ_s \Vert^2 ds
   \le C \left( \int_t^T \esp \vert \dX_s \vert^2 ds + A_{l,m} \right).
\label{estimdz}
\end{equation}
Gronwall's lemma gives $\esp \vert \dX_t \vert^2 \le C A_{l,m}$, where
$C$ is again independent of $l$ and $m$.
Moreover, using \eqref{lip}, $f_l$ converges uniformly to $f$ on
$\R^d \times \omb \times B(0,r)$ for
any $r>0$ (with $B(0,r)= \{ z \in \R^{n d_W} : \Vert z \Vert < r \}$) and
$X^l$, $X^m$ are bounded, so
$A_{l,m}$ tends to zero when $l,m$ tend to infinity; therefore $(X^l)_l$
converges to a process $X$ in $L^2(\Omega \times [0;T])$.

Using \eqref{estimdz} again, we get
$$\esp \int_0^T \Vert \dZ_s \Vert^2 ds \le C A_{l,m};$$
hence the sequence of processes $(Z^l)$ has also a limit in
$L^2(\Omega \times [0;T])$; let $Z$ denote this limit process.

The pair $(X,Z)$ solves BSDE $(M+D)$ and $X$ is $\omb$-valued; the
proof is just an adaptation of the second step in the proof of
Proposition \ref{passlim}.
This remark completes the proof.
\end{demo}

{\bf Remark :} As a consequence, a.s.,
$\Vert Z_t \Vert \le \frac{1}{\ep}$ for any $t$.

According to Subsection \ref{par4.1}, this result can be extended to
every $\omb$-valued and ${\cal F}_T$-measurable terminal variable $U$.
Then Theorem \ref{existence2} of existence and uniqueness of a solution
follows.

Note that uniqueness and existence hold in particular on any regular
geodesic ball (or geodesic ball if the sectional curvatures are
nonpositive).

\section{Applications and related PDEs}
\label{par5}

\subsection{The martingale case}
\label{par5.1}
The drift $f=0$ verifies hypothesis $(H)$. Hence in this case the results
of this paper apply to the martingale case. As already underlined, any
regular geodesic ball verifies the condition of Theorem
\ref{existence2}; so we recover the well-known results of existence and
uniqueness of a martingale with prescribed terminal value in such domains
(see \cite{kend90}).
These results hold in nonpositive curvatures, they will be achieved in
positive curvatures elsewhere.

\subsection{The one-dimensional case}
\label{par5.2}
The nonpositive curvature case gives the existence and uniqueness of a
solution to the one-dimensional BSDE
$$
(E)_1
\left\{
      \begin{array}{l}
           d X_t = Z_t d W_t - \Gamma(X_t) Z_t^2 + f(B_t^y,X_t,Z_t) dt \\
           X_T=U
      \end{array}
      \right.
$$
for a bounded terminal condition $U$, a drift $f$ satisfying
\eqref{upperboundf}, \eqref{lip2}
and any smooth function $\Gamma$ defined on $\R$.

Note that a change of coordinates (in fact a reparametrization of the
one-dimensional manifold by arclength) reduces equation $(E)_1$ to
$$
\left\{
      \begin{array}{l}
           d X_t = Z_t d W_t + f(B_t^y,X_t,Z_t) dt  \\
           X_T=U;
      \end{array}
      \right.
$$
moreover, it is a very particular case of the results of Kobylansky
in \cite{kobyl00}.

One can ask whether such results can be extended to higher dimensions.
In fact, the original problem is geometric and to deal with general BSDEs,
we would start with smooth functions $(\Gamma_{ij}^k)$ and should give
conditions in order to interpret these functions as the Christoffel
symbols of a given Levi-Civita connection. This problem is out of the
scope of this paper.

\subsection{Case of a random terminal time}
\label{par5.3}
In this paragraph, we will only sketch the proofs.

We are interested in the following equation
$$(M+D)_\tau
\left\{
      \begin{array}{l}
          d X_t = Z_t d W_t + \left( - \half \Gamma_{jk}(X_t)
            ([Z_t]^k \vert [Z_t]^j) + f(B_t^y,X_t,Z_t) \right) d t  \\
	  X_\tau=U^\tau                                         \\
      \end{array}
    \right.   $$
where $\tau$ is a stopping time with respect to the filtration used and
$U^\tau$ is a $\omb$-valued, ${\cal F}_\tau$-measurable random variable.
It is the counterpart of equation $(M+D)$ on the random interval
$[0; \tau]$.

Let us first consider the case of a bounded stopping time $\tau$, i.e.
$\tau \le T$ where $T$ is a deterministic constant.
We have the following result :

\begin{stheorem}
\label{existence3}
We consider BSDE $(M+D)_\tau$ with $\omb= \{ \chi \le c \}$ and
$\tau \le T$ a.s. If $f$
verifies conditions \eqref{lip}, \eqref{upperboundf} and $(H)$, and
if $\chi$ is strictly convex (i.e. $\Hess \chi$ is positive definite),
then this BSDE has a unique solution $(X,Z)$, with $X \in \omb$, in
the same two cases as in Theorem \ref{existence2}.
\end{stheorem}

\begin{demo}
First remark that the uniqueness part goes the same as in the
deterministic case; for the existence part, it can be completed in
several steps~:

{\it First Step : }
We work in local coordinates introduced in Subsection \ref{par4.4}.
Let $c_1>c$ and put $\omb_1= \{ \chi \le c_1 \}$; suppose that $c_1$
is such that $\omb \subset \omb_1 \subset O$ and that $\chi$ is yet
strictly convex on $\omb_1$. Now let $\phi$ be a cut-off function with
$\phi = 1$ on $\omb$ and $\phi = 0$ outside $\omb_1$.
For any nonzero integer $l$, we define a new drift by
$f_l(b,x,z):=\phi(x) (f(b,x,z) + (1/l) x)$; note
that $f_l$
verifies hypothesis $(H)$ with respect to $\omb_1$ and $(H_s)$ with
respect to $\omb$.

Then solve path by path on $[\tau;T]$ the following differential equation
\begin{equation}
\left\{
      \begin{array}{l}
          d X_t^l = f_l(B_t^y,X_t^l, 0) d t  \\
	  X_\tau^l=U^\tau
      \end{array}
    \right.
\label{eqtaut}
\end{equation}
and set $U^l= U^\tau + \int_\tau^T f_l(B_t^y,X_t^l,0) dt$. Since $f_l$
vanishes outside $\omb_1$, $U^l$ is in $\omb_1$; besides, $U^l$ is
${\cal F}_T$-measurable.

{\it Second Step : }
Considering the random variable $U^l$ and the drift $f_l$ introduced in
the first step, we solve on $[0;T]$
equation $(M+D)$ with drift $f_l$ and terminal value
$X_T^l=U^l$. The hypothesis of Theorem \ref{existence2} are satisfied
considering the set $\omb_1$ instead of $\omb$. So this BSDE has a
solution $(X^l,Z^l)$ with $X^l \in \omb_1$.

We condition by ${\cal F}_\tau$ and consider the above solution $(X^l,Z^l)$
on the random time interval $[\tau; T]$. It is a solution on this interval
of BSDE $(M+D)$ with drift $f_l$ and terminal value $U^l$. The uniqueness
property for such equations implies that $X^l$ is the solution
of equation \eqref{eqtaut}, i.e.
$$\forall t \in [\tau;T],\ X_t^l = U^l - \int_t^T f_l(B_s^y,X_s^l,0) ds$$
and $Z_t^l = 0$ for $\tau \le t \le T$. In particular,
$$X_\tau^l = U^l - \int_\tau^T f_l(B_s^y,X_s^l,0) ds = U^\tau.$$

We now show that actually, $(X_t^l)_{0 \le t \le \tau}$ remains in $\omb$,
and not only in $\omb_1$.
For this purpose, we want to construct as in Subsection \ref{par4.4} a
submartingale which is written as $(e^{\la t} H(X_t^l))$, where now
$H=h \circ \chi$, since $X^l$ is $\omb_1$-valued.
Recall from \eqref{pos3} that the keypoint is to show for $x \in \omb_1$
the nonnegativity of
$$
h'(\chi(x)) \left(
\half
  \sum_{i=1}^{d_W} {}^t [{}^t z]^i \Hess \chi (x) [{}^t z]^i
    +  D \chi(x) f_l(b,x,z) \right) + \lambda H(x).
$$
In fact, as $\chi$ is strictly convex on the compact domain $\omb_1$,
we have $\Hess \chi \ge \al Id$ (in the sense of matrices) for $\al>0$
and it turns out that it suffices to prove the nonnegativity of

\begin{equation}
h'(\chi(x)) \left( \half
  \al \Vert z \Vert^2 +  D \chi(x) f_l(b,x,z) \right) + \lambda H(x).
\label{pos4}
\end{equation}

From \eqref{lip} and \eqref{upperboundf} we deduce
$$ \vert D \chi(x) f_l(b,x,z) \vert \le C (1+ \Vert z \Vert) \le
   \half \al \Vert z \Vert^2,$$
the last inequality holding for $\Vert z \Vert$ large enough, say
$\Vert z \Vert \ge A$. Obviously in this case, \eqref{pos4} is
nonnegative. \\
Now suppose that $\Vert z \Vert \le A$. If
$x \in \overline{B(0,1)}=\omb$ then $h'=0$ and the required result
holds. Otherwise we write for $x \in \omb_1 \setminus \omb$ and
$x_0 \in \partial \omb$ (i.e. in our local coordinates
$\vert x \vert \ge 1$ and $\vert x_0 \vert =1$)~:
$$ \vert D \chi(x) f_l(b,x,z)- D \chi(x_0) f_l(b,x_0,z) \vert \le C
    \vert x-y \vert$$
for a constant $C$. But the hypothesis $(H_s)$ for $f_l$ writes
$$ D \chi(x_0) f_l(b,x_0,z) \ge \zeta >0.$$
Then we distinguish two cases ($x$ near $1$ or "far" from $1$) and get
the nonnegativity of \eqref{pos4} in both situations. This can be done
by using similar arguments to those displayed at the end of the Third
Step in the proof of Proposition \ref{prop4.1}; in particular, a $\la$
large enough is needed.

{\it Third Step : }
The first two steps give the existence of processes $(X^l,Z^l)$ (with
$X^l \in\omb$) solving equation $(M+D)_\tau$ associated to the drift
$f^l$ and terminal value $U^\tau$. But $f^l$ converges to $f$ uniformly
on $\omb$ so, passing through the limit as in Subsection \ref{par4.6},
we get a pair of processes $(X_t,Z_t)_{0 \le t \le \tau}$, with
$X \in \omb$ and solving the initial equation $(M+D)_\tau$ with drift
$f$ and terminal value $U^\tau$.
This completes the proof.
\end{demo}

\bigskip
We consider again a stopping time $\tau$ and the corresponding equation
$(M+D)_\tau$; now, we only suppose that $\tau$ is finite a.s. and
verifies the exponential integrability condition \eqref{integst}.
Examples of such stopping times are exit times of uniformly elliptic
diffusions from bounded domains in Euclidean spaces.

In this case, we need to add restrictions on the drift $f$; indeed,
the main thrust in the proof of uniqueness and existence is the
construction of a submartingale on the product manifold
$(S_t)_t = (\exp(A_t) \Psi(\tilde X_t))_t$
with $\mu = 0$ if $f$ does not depend on $z$. To extend this approach
to a random (non necessarily bounded) interval, we have to keep the
integrability of $S_\tau$.
An accurate examination of the method to obtain the submartingale (in
particular inequalities \eqref{majdpsi} and \eqref{integ}) shows that
this integrability holds for "small" drifts; more precisely there is
a constant $h$ with $0<h<\rho$ such that, under the following
condition on the constants in \eqref{lip} and \eqref{upperboundf}
\begin{equation}
L < h, \ \ L_2 < h,
\label{smalldrift}
\end{equation}
the integrability required holds, so $(S_t)_{0 \le t \le \tau}$ is a
true submartingale.

{\bf Remarks : }
1- Such a condition guarantees in particular that we have
$$ \esp \int_0^\tau \vert f(B_s^y,X_s,Z_s) \vert ds < \infty.$$

2- This condition is rather natural; actually, it is very similar to
conditions yet introduced for BSDEs with Lipschitz coefficients and
random terminal time~: see condition (24) and Propositions 3.2 and
3.3 in \cite{darlpard97}, or (2.6) and the condition before in
\cite{peng91}.

3- {\it A priori}, the process $(Z_t)_t$ verifies the integrability
condition
$$\esp \left( \int_0^\tau \Vert Z_s \Vert^2 ds \right) < \infty;$$
in fact, it results from the existence part that in any case (i.e.
$f$ depending or not on $z$), $(Z_t)_{0 \le t \le \tau}$ belongs to
$(\cal E_\al)$
(see Definition \ref{d1}) for $\al$ small enough, which is a
stronger property. In particular, we get
$$ \forall \theta < \rho, \ \ \esp \int_0^\tau e^{\theta s}
    \Vert Z_s \Vert^2 ds < \infty;$$
this condition is usual for BSDEs with random terminal time (see again
Propositions 3.2 and 3.3 in \cite{darlpard97}, or Theorem 2.2 in
\cite{peng91}).

Once we have constructed the submartingale as on a deterministic
interval, uniqueness is straightforward. Let us indicate how existence
can be deduced.

We are given a $\omb$-valued and ${\cal F}_\tau$-measurable variable
$U^\tau$.
As in the proof of Theorem \ref{existence3}, we consider again local
coordinates introduced in Subsection \ref{par4.4}, $\omb_1$ such that
$\omb \subset \omb_1 \subset O$ and a cut-off function $\phi$ with
$\phi = 1$ on $\omb$ and $\phi = 0$ outside $\omb_1$. We put
$f_1(b,x,z)=\phi(x)f(b,x,z)$. \\
The first step here is to solve on $[0;\tau]$ a BSDE whose terminal
value is near $U^\tau$~: \\
on $[0;\tau \wedge n]$, using Theorem \ref{existence3}, we solve
equation $(M+D)$ with drift $f_1$ and terminal value at time
$\tau \wedge n$, $\esp [U^\tau \vert {\cal F}_n]$; let
$(X_t^n, Z_t^n)_{0 \le t \le \tau \wedge n}$ denote the solution; \\
on $[\tau \wedge n ; \tau]$, we put $Z_t=0$ and solve
$$
\left\{
      \begin{array}{l}
          d X_t^n = f_1(B_t^y,X_t^n, 0) d t  \\
	  X_{\tau \wedge n}^n=\esp [U^\tau \vert {\cal F}_n].
      \end{array}
    \right.
$$
Then it is easily seen, since $\esp [U^\tau \vert {\cal F}_n]$ is
${\cal F}_{\tau \wedge n}$-measurable, that
$(X_t^n, Z_t^n)_{0 \le t \le \tau}$ is a solution to BSDE $(M+D)_\tau$
with terminal value $U^{\tau,n}$, where
$$U^{\tau,n} = X_{\tau \wedge n}^n + \int_{\tau \wedge n}^\tau
   f_1(B_s^y,X_s^n,0) ds.$$
The second step is to show that when $n$ tends to infinity, we get
the solution of BSDE $(M+D)_\tau$ with terminal value $U^\tau$. \\
We have that $U^{\tau,n}$ tends to $U^\tau$ in $L^2(\Omega)$; indeed,
$$\esp \vert U^\tau-U^{\tau,n} \vert ^2 =
    \esp \left( 1_{n \le \tau} \left\vert U^\tau
    - \esp [U^\tau \vert {\cal F}_n] - \int_{\tau \wedge n}^\tau
    f_1(B_s^y,X_s^n,0) ds \right\vert ^2 \right)$$
and the last expectation tends to zero as $n$ tends to infinity; this is a
consequence of dominated convergence, using the exponential integrability
condition of $\tau$. \\
Then we apply the results of
Subsection \ref{par4.1} to the random interval $[0;\tau]$; indeed, an
accurate examination shows that these results rely essentially on
the exponential integrability condition (according to the uniqueness part)
$$\esp \left( e^{\al \int_0^\tau \left( \Vert Z_s^n \Vert^2
    + \Vert Z_s^m \Vert^2 \right) ds} \right) \le C < \infty$$
for a constant $C$ independent of $m,n$.
At the end, we get the existence of a pair of processes
$(X_t, Z_t)_{0 \le t \le \tau}$ solution of BSDE $(M+D)_\tau$ with
drift $f_1$ and terminal value $U^\tau$. But, since $U^\tau$ is
$\omb$-valued and $\omb=\overline{B(0,1)}$, for each $n$ the process
$(X_t^n)_{0 \le t \le \tau \wedge n}$ remains in $\omb$ by Theorem
\ref{existence3}; thus the whole process $(X_t)_{0 \le t \le \tau}$
remains in $\omb$ and this completes the existence part.

\medskip
As a consequence, we can state the following result
\begin{stheorem}
\label{existence4}
We consider BSDE $(M+D)_\tau$ with $\tau$ a stopping time verifying the
integrability condition \eqref{integst}; the function $\chi$ used to
define the domain $\omb$ is supposed as usual to be strictly convex.
Then if $f$ verifies conditions \eqref{lip}, \eqref{upperboundf}, $(H)$
and moreover is "small" (i.e. verifies condition \eqref{smalldrift}
above), this BSDE has a unique solution
$(X,Z)$, in the same cases as in Theorem \ref{existence2}.
\end{stheorem}

\subsection{Application to nonlinear elliptic PDEs}
\label{par5.4}
In this paragraph, we make precise the Dirichlet problem that we briefly
discussed in the introduction.

Suppose $(N,g)$ is a Riemannian manifold, and $B^x$ a Brownian motion
on $(N,g)$ (started at $x$ at time $0$). Alternatively, think of $B^x$
as the diffusion process on $\R^d$, defined by \eqref{sde}; in this case,
$$\forall i,j=1, \ldots, d, \ \ \sum_{l=1}^d \si_{il} \si_{jl} = g^{ij},$$
the inverse metric tensor, and
$$\forall i=1, \ldots, d, \ \
      b^i + \sum_{k,l=1}^{d_W} g^{kl} \Gamma^i_{kl}=0.$$
Let $\overline M_1$ be a compact submanifold of $N$, with boundary
$\partial M_1$ and interior $M_1$. For $x \in \overline M_1$, we call
$\zeta$ the first time $B^x$
hits the boundary; we assume that $\zeta$ verifies an integrability
condition like \eqref{integst}.
Given a regular mapping
$$\overline \phi : \partial M_1 \rightarrow \omb \subset M,$$
we wish to find a mapping $\phi : \overline M_1 \rightarrow \omb$
which solves the following Dirichlet problem
$$(D)
\left\{
\begin{array}{ccc}
{\cal L}_M \phi(x) - f(x, \phi(x), \nabla \phi(x) \si(x))=0
& , & x \in M_1                                        \\
\phi(x) = \overline \phi (x)
& , & x \in \partial M_1
\end{array}
\right.
$$
where ${\cal L}_M \phi$ is the tension field of the mapping $\phi$ (see
\cite{eellsamp}, or for a probabilistic point of view the introduction
of \cite{mpl}).

We recall from the introduction that, in coordinates $(x^i)$ on $M$ and
$(y^\alpha)$ on $M_1$, the equation ${\cal L}_M \phi = 0$ characterizes
harmonic mappings, and is written
$$ \forall i, \ \ \Delta_{M_1} \phi^i + g^{\al \beta} \Gamma^i_{jk}(\phi)
      D_\alpha \phi^j D_\beta \phi^k =0.$$
Using the same Wiener process $W$ with which we constructed $B^x$, we can
solve according to Theorem \ref{existence4} the BSDE $(M+D)_\zeta$ with
terminal value $\overline \phi (B^x_\zeta)$. Let
$(X_t^x, Z_t^x)_{0 \le t \le \zeta}$ be the unique solution and put
$\phi(x):=X_0^x$. Then
under sufficient regularity on $\phi$, it is not difficult to verify that
$\phi$ is a solution to the Dirichlet problem $(D)$. Note that when
$f \equiv 0$ (i.e. in the martingale case), Kendall (\cite{kend94}) has
proved regularity results on $\phi$ using almost only probability theory,
so that $\phi$ is a strong solution of the equation
${\cal L}_M \phi \equiv 0$ (i.e. a harmonic mapping).

When $f(b,x,z)=f(b,x)$ and is written as
$f(b,x)  = D_2 G(b,x)$ (the differential of $G$ with respect to the
second variable), the elliptic nonlinear PDE in the Dirichlet
problem $(D)$ is associated with a variational problem; more precisely,
solutions of this equation are critical points of the functional
$$ {\cal F} (u) = \half  \int_{M_1} \Vert {\rm grad} u(b) \Vert ^2
      d vol(b) + \int_{M_1} G (b,u(b)) d vol(b) $$
and the elliptic PDE in equation $(D)$ is the Euler-Lagrange equation
associated.

\subsection{Application to nonlinear parabolic PDEs}
\label{par5.5}
We conclude this part by studying the time-dependent equation associated
with the stationary equation described in the Dirichlet problem $(D)$
above. More precisely, we are interested in the following equation, for
mappings $u : [0;T] \times N \rightarrow \omb \subset M$~:
$$
\left\{
\begin{array}{rcl}
\frac{\partial u}{\partial t}
& = & {\cal L}_M u - f(x, u, \nabla u \si)  \\
u_{|t=0}
& = & F
\end{array}
\right.
$$
where $F$ is sufficiently regular and has range $\omb$. In local
coordinates, this equation becomes
$$
\left\{
\begin{array}{rclcl}
\frac{\partial u}{\partial t}(t,x)
& = & \half {\cal L} u(t,x)
& + & \half \Gamma_{jk}(u(t,x))
       ([(\nabla_x u \si)(t,x)]^k \vert [(\nabla_x u \si)(t,x)]^j)  \\
&   &
& - & f(x,u(t,x), (\nabla_x u \si)(t,x))    \\
u(0,x)
& = & F(x)
\end{array}
\right.
$$
in the case of ${\cal L}$ being the Laplace-Beltrami operator on $N$. This
is equation \eqref{eqchaleur}.
As a by-product of Section \ref{par4}, we have the existence and
uniqueness of a regular solution to this system of quasilinear parabolic
PDEs; it is based on the boundedness of $\nabla_x u$, proved in
Subsection \ref{par4.3}.

\bibliography{BSDE1}
\bibliographystyle{abbrv}

\end{document}